\documentclass[12pt]{amsart}      %{article} was 12pt latex e
\usepackage{amssymb}
\usepackage{eucal}
\usepackage{amsmath}
\usepackage{amscd}
\usepackage{multicol}
\usepackage[all]{xy}           %xypic macro for latex2.09
\usepackage{graphicx}
\usepackage{xspace}
\usepackage{axodraw}
\begin{document}  %for latex 2.09

\newcommand{\nc}{\newcommand}
\newcommand{\delete}[1]{}
\nc{\dfootnote}[1]{{}}          %{{}}
\nc{\ffootnote}[1]{\dfootnote{#1}}
%\nc{\mfootnote}[1]{{}}        % Use this to suppress footnotes
\nc{\mfootnote}[1]{\footnote{#1}} % Use this to show footnotes
%\nc{\ofootnote}[1]{{}}        % Use this to suppress footnotes
\nc{\ofootnote}[1]{\footnote{\tiny Older version: #1}} % Use this to show footnotes

\nc{\mlabel}[1]{\label{#1}}  % Use this to suppress names
\nc{\mcite}[1]{\cite{#1}}  % Use this to suppress names
\nc{\mref}[1]{\ref{#1}}  % Use this to suppress names

\delete{
\nc{\mlabel}[1]{\label{#1}  % Use the next two lines to show names
{\hfill \hspace{1cm}{\bf{{\ }\hfill(#1)}}}}
\nc{\mcite}[1]{\cite{#1}{{\bf{{\ }(#1)}}}}  % Use this lines to show names
\nc{\mref}[1]{\ref{#1}{{\bf{{\ }(#1)}}}}  % Use this lines to show names
}
\nc{\mbibitem}[1]{\bibitem{#1}} % Use this to show number
%\nc{\mbibitem}[1]{\bibitem[\bf #1]{#1}} % Use this to show name
\nc{\mkeep}[1]{\marginpar{{\bf #1}}} % Use this to show marginpar
%\nc{\mkeep}[1]{{}}      % Use this to suppress marginpar

%%%%%%%%%%%%%%%%%%%%%%%% Statements
\newtheorem{theorem}{Theorem}[section]
\newtheorem{prop}[theorem]{Proposition}
\newtheorem{defn}[theorem]{Definition}
\newtheorem{lemma}[theorem]{Lemma}
\newtheorem{coro}[theorem]{Corollary}
\newtheorem{prop-def}{Proposition-Definition}[section]
\newtheorem{claim}{Claim}[section]
\newtheorem{remark}[theorem]{Remark}
\newtheorem{propprop}{Proposed Proposition}[section]
\newtheorem{conjecture}{Conjecture}
\newtheorem{exam}{Example}[section]
\newtheorem{assumption}{Assumption}
\newtheorem{condition}[theorem]{Assumption}

\renewcommand{\labelenumi}{{\rm(\alph{enumi})}}
\renewcommand{\theenumi}{\alph{enumi}}

\nc{\tred}[1]{\textcolor{red}{#1}}
\nc{\tblue}[1]{\textcolor{blue}{#1}}
\nc{\tgreen}[1]{\textcolor{green}{#1}}

%%%%%%%%%%%%%%%%%%%%%%% symbols
\nc{\adec}{\check{;}}
\nc{\dftimes}{\widetilde{\otimes}} \nc{\dfl}{\succ}
\nc{\dfr}{\prec} \nc{\dfc}{\circ} \nc{\dfb}{\bullet}
\nc{\dft}{\star} \nc{\dfcf}{{\mathbf k}} \nc{\spr}{\cdot}
\nc{\disp}[1]{\displaystyle{#1}}
\nc{\bin}[2]{ (_{\stackrel{\scs{#1}}{\scs{#2}}})}  %binomial coeff
\nc{\binc}[2]{ \left (\!\! \begin{array}{c} \scs{#1}\\
    \scs{#2} \end{array}\!\! \right )}  %binomial coeff
\nc{\bincc}[2]{  \left ( {\scs{#1} \atop
    \vspace{-.5cm}\scs{#2}} \right )}  %binomial coeff
\nc{\sarray}[2]{\begin{array}{c}#1 \vspace{.1cm}\\ \hline
    \vspace{-.35cm} \\ #2 \end{array}}
\nc{\bs}{\bar{S}} \nc{\dcup}{\stackrel{\bullet}{\cup}}
\nc{\dbigcup}{\stackrel{\bullet}{\bigcup}} \nc{\etree}{\big |}
\nc{\la}{\longrightarrow} \nc{\fe}{\'{e}} \nc{\rar}{\rightarrow}
\nc{\dar}{\downarrow} \nc{\dap}[1]{\downarrow
\rlap{$\scriptstyle{#1}$}} \nc{\uap}[1]{\uparrow
\rlap{$\scriptstyle{#1}$}} \nc{\defeq}{\stackrel{\rm def}{=}}
\nc{\dis}[1]{\displaystyle{#1}} \nc{\dotcup}{\,
\displaystyle{\bigcup^\bullet}\ } \nc{\sdotcup}{\tiny{
\displaystyle{\bigcup^\bullet}\ }} \nc{\hcm}{\ \hat{,}\ }
\nc{\hcirc}{\hat{\circ}} \nc{\hts}{\hat{\shpr}}
\nc{\lts}{\stackrel{\leftarrow}{\shpr}}
\nc{\rts}{\stackrel{\rightarrow}{\shpr}} \nc{\lleft}{[}
\nc{\lright}{]} \nc{\uni}[1]{\tilde{#1}} \nc{\wor}[1]{\check{#1}}
\nc{\free}[1]{\bar{#1}} \nc{\den}[1]{\check{#1}} \nc{\lrpa}{\wr}
\nc{\curlyl}{\left \{ \begin{array}{c} {} \\ {} \end{array}
    \right .  \!\!\!\!\!\!\!}
\nc{\curlyr}{ \!\!\!\!\!\!\!
    \left . \begin{array}{c} {} \\ {} \end{array}
    \right \} }
\nc{\leaf}{\ell}       % number of leafs
\nc{\longmid}{\left | \begin{array}{c} {} \\ {} \end{array}
    \right . \!\!\!\!\!\!\!}
\nc{\ot}{\otimes} \nc{\sot}{{\scriptstyle{\ot}}}
\nc{\otm}{\overline{\ot}}
\nc{\ora}[1]{\stackrel{#1}{\rar}}
\nc{\ola}[1]{\stackrel{#1}{\la}}%${\Bbb Z}$
\nc{\scs}[1]{\scriptstyle{#1}} \nc{\mrm}[1]{{\rm #1}}
\nc{\margin}[1]{\marginpar{\rm #1}}   %{\rm #1}}
\nc{\dirlim}{\displaystyle{\lim_{\longrightarrow}}\,}
\nc{\invlim}{\displaystyle{\lim_{\longleftarrow}}\,}
\nc{\mvp}{\vspace{0.5cm}} \nc{\svp}{\vspace{2cm}}
\nc{\vp}{\vspace{8cm}} \nc{\proofbegin}{\noindent{\bf Proof: }}
%\nc{\proofbegin}{\begin{proof}} % AMS command
\nc{\proofend}{$\blacksquare$ \vspace{0.5cm}}
%\nc{\proofend}{\end{proof}} %AMS command
%\nc{\intg}[1]{\lceil{#1}\rceil}  %old free int ring
%\nc{\sha}{\scs{\mbox{\cyr X}}} %used to be \cyr
\nc{\sha}{{\mbox{\cyr X}}}  %used to be \cyr
\nc{\ncsha}{{\mbox{\cyr X}^{\mathrm NC}}} \nc{\ncshao}{{\mbox{\cyr
X}^{\mathrm NC,\,0}}}
\nc{\shpr}{\diamond}    %Shuffle product
\nc{\shprm}{\overline{\diamond}}    %Shuffle product
\nc{\shpro}{\diamond^0}    %Shuffle product
\nc{\shprr}{\diamond^r}     %product on controlled trees
\nc{\shpra}{\overline{\diamond}^r}
\nc{\shpru}{\check{\diamond}} \nc{\catpr}{\diamond_l}
\nc{\rcatpr}{\diamond_r} \nc{\lapr}{\diamond_a}
\nc{\sqcupm}{\ot}
\nc{\lepr}{\diamond_e} \nc{\vep}{\varepsilon} \nc{\labs}{\mid\!}
\nc{\rabs}{\!\mid} \nc{\hsha}{\widehat{\sha}}
\nc{\lsha}{\stackrel{\leftarrow}{\sha}}
\nc{\rsha}{\stackrel{\rightarrow}{\sha}} \nc{\lc}{\lfloor}
\nc{\rc}{\rfloor} \nc{\sqmon}[1]{\langle #1\rangle}
\nc{\forest}{\calf} \nc{\ass}[1]{\alpha({#1})}
\nc{\altx}{\Lambda_X} \nc{\vecT}{\vec{T}} \nc{\onetree}{\bullet}
\nc{\Ao}{\check{A}}
\nc{\seta}{\underline{\Ao}}
\nc{\deltaa}{\overline{\delta}}
\nc{\trho}{\tilde{\rho}}
\nc{\tpow}[2]{{#2}^{\ot #1}}

%%%%%%%%%%%%%%%%%%%%% roman fonts, in alphabetic order
\nc{\mmbox}[1]{\mbox{\ #1\ }} \nc{\ann}{\mrm{ann}}
\nc{\Aut}{\mrm{Aut}} \nc{\can}{\mrm{can}} \nc{\colim}{\mrm{colim}}
\nc{\Cont}{\mrm{Cont}} \nc{\rchar}{\mrm{char}}
\nc{\cok}{\mrm{coker}} \nc{\dtf}{{R-{\rm tf}}} \nc{\dtor}{{R-{\rm
tor}}}
\renewcommand{\det}{\mrm{det}}
\nc{\depth}{{\mrm d}}
\nc{\Div}{{\mrm Div}} \nc{\End}{\mrm{End}} \nc{\Ext}{\mrm{Ext}}
\nc{\Fil}{\mrm{Fil}} \nc{\Frob}{\mrm{Frob}} \nc{\Gal}{\mrm{Gal}}
\nc{\GL}{\mrm{GL}} \nc{\Hom}{\mrm{Hom}} \nc{\hsr}{\mrm{H}}
\nc{\hpol}{\mrm{HP}} \nc{\id}{\mrm{id}} \nc{\im}{\mrm{im}}
\nc{\incl}{\mrm{incl}} \nc{\length}{\mrm{length}}
\nc{\LR}{\mrm{LR}} \nc{\mchar}{\rm char} \nc{\NC}{\mrm{NC}}
\nc{\mpart}{\mrm{part}} \nc{\ql}{{\QQ_\ell}} \nc{\qp}{{\QQ_p}}
\nc{\rank}{\mrm{rank}} \nc{\rba}{\rm{RBA }} \nc{\rbas}{\rm{RBAs }}
\nc{\rbw}{\rm{RBW }} \nc{\rbws}{\rm{RBWs }} \nc{\rcot}{\mrm{cot}}
\nc{\rest}{\rm{controlled}\xspace}
\nc{\rdef}{\mrm{def}} \nc{\rdiv}{{\rm div}} \nc{\rtf}{{\rm tf}}
\nc{\rtor}{{\rm tor}} \nc{\res}{\mrm{res}} \nc{\SL}{\mrm{SL}}
\nc{\Spec}{\mrm{Spec}} \nc{\tor}{\mrm{tor}} \nc{\Tr}{\mrm{Tr}}
\nc{\mtr}{\mrm{sk}}

%%%%%%%%%%%%%%%%%% bold face
\nc{\ab}{\mathbf{Ab}} \nc{\Alg}{\mathbf{Alg}}
\nc{\Algo}{\mathbf{Alg}^0} \nc{\Bax}{\mathbf{Bax}}
\nc{\Baxo}{\mathbf{Bax}^0} \nc{\RB}{\mathbf{RB}}
\nc{\RBo}{\mathbf{RB}^0} \nc{\BRB}{\mathbf{RB}}
\nc{\Dend}{\mathbf{DD}} \nc{\bfk}{{\bf k}} \nc{\bfone}{{\bf 1}}
\nc{\base}[1]{{a_{#1}}} \nc{\detail}{\marginpar{\bf More detail}
    \noindent{\bf Need more detail!}
    \svp}
\nc{\Diff}{\mathbf{Diff}} \nc{\gap}{\marginpar{\bf
Incomplete}\noindent{\bf Incomplete!!}
    \svp}
\nc{\FMod}{\mathbf{FMod}} \nc{\mset}{\mathbf{MSet}}
\nc{\rb}{\mathrm{RB}} \nc{\Int}{\mathbf{Int}}
\nc{\Mon}{\mathbf{Mon}}
%\nc{\remark}{\noindent{\bf Remark: }}
\nc{\remarks}{\noindent{\bf Remarks: }} \nc{\Rep}{\mathbf{Rep}}
\nc{\Rings}{\mathbf{Rings}} \nc{\Sets}{\mathbf{Sets}}
\nc{\DT}{\mathbf{DT}}

%%%%%%%%%%%%%%%%%%%Bbb fonts
\nc{\BA}{{\mathbb A}} \nc{\CC}{{\mathbb C}} \nc{\DD}{{\mathbb D}}
\nc{\EE}{{\mathbb E}} \nc{\FF}{{\mathbb F}} \nc{\GG}{{\mathbb G}}
\nc{\HH}{{\mathbb H}} \nc{\LL}{{\mathbb L}} \nc{\NN}{{\mathbb N}}
\nc{\QQ}{{\mathbb Q}} \nc{\RR}{{\mathbb R}} \nc{\TT}{{\mathbb T}}
\nc{\VV}{{\mathbb V}} \nc{\ZZ}{{\mathbb Z}}

%%%%%%%%%%%%%%%%%%% cal fonts

\nc{\cala}{{\mathcal A}} \nc{\calc}{{\mathcal C}}
\nc{\cald}{{\mathcal D}} \nc{\cale}{{\mathcal E}}
\nc{\calf}{{\mathcal F}} \nc{\calfr}{{{\mathcal F}^{\,r}}}
\nc{\calfo}{{\mathcal F}^0} \nc{\calfro}{{\mathcal F}^{\,r,0}}
\nc{\oF}{\overline{F}}  \nc{\calg}{{\mathcal G}}
\nc{\calh}{{\mathcal H}} \nc{\cali}{{\mathcal I}}
\nc{\calj}{{\mathcal J}} \nc{\call}{{\mathcal L}}
\nc{\calm}{{\mathcal M}} \nc{\caln}{{\mathcal N}}
\nc{\calo}{{\mathcal O}} \nc{\calp}{{\mathcal P}}
\nc{\calr}{{\mathcal R}} \nc{\calt}{{\mathcal T}}
\nc{\caltr}{{\mathcal T}^{\,r}}
\nc{\calu}{{\mathcal U}} \nc{\calv}{{\mathcal V}}
\nc{\calw}{{\mathcal W}} \nc{\calx}{{\mathcal X}}
\nc{\CA}{\mathcal{A}}

%%%%%%%%%%%%%%%%%%  frak fonts
\nc{\fraka}{{\mathfrak a}} \nc{\frakB}{{\mathfrak B}}
\nc{\frakb}{{\mathfrak b}} \nc{\frakd}{{\mathfrak d}}
\nc{\oD}{\overline{D}}
\nc{\frakF}{{\mathfrak F}} \nc{\frakg}{{\mathfrak g}}
\nc{\frakm}{{\mathfrak m}} \nc{\frakM}{{\mathfrak M}}
\nc{\frakMo}{{\mathfrak M}^0} \nc{\frakp}{{\mathfrak p}}
\nc{\frakS}{{\mathfrak S}} \nc{\frakSo}{{\mathfrak S}^0}
\nc{\fraks}{{\mathfrak s}} \nc{\os}{\overline{\fraks}}
\nc{\frakT}{{\mathfrak T}}
\nc{\oT}{\overline{T}}
%\nc{\frakx}{{\mathfrak x}}
\nc{\frakX}{{\mathfrak X}} \nc{\frakXo}{{\mathfrak X}^0}
\nc{\frakx}{{\mathbf x}}
%\nc{\frakTxo}{{\frakTx}^0}
\nc{\frakTx}{\frakT}      %All rooted trees, correspond to \ncsha(X)
\nc{\frakTa}{\frakT^a}        % rooted trees for \ncsha(A)
\nc{\frakTxo}{\frakTx^0}   % rooted trees for \ncshao(X)
\nc{\caltao}{\calt^{a,0}}   % rooted trees for \ncshao(A)
\nc{\ox}{\overline{\frakx}} \nc{\fraky}{{\mathfrak y}}
\nc{\frakz}{{\mathfrak z}} \nc{\oX}{\overline{X}}

\font\cyr=wncyr10

\nc{\redtext}[1]{\textcolor{red}{#1}}

\nc{\arbreA}{%\kern-0.4ex
\setlength{\unitlength}{.7pt}
\begin{picture}(60,40)(0,0)
\put(30,0){\line(0,1){10}}
\put(30,10){\line(-1,1){30}}
\put(30,10){\line(1,1){30}}
\end{picture}}%\kern 0.4ex}

\nc{\arbreAd}{%\kern-0.4ex
\setlength{\unitlength}{.7pt}
\begin{picture}(60,40)(0,0)
\put(30,0){\line(0,1){10}}
\put(30,10){\line(-1,1){30}}
\put(30,10){\line(1,1){30}}
\put(35,3){\small{x}}
\end{picture}}%\kern 0.4ex}

\def\arbreB{%\kern-0.4ex
\setlength{\unitlength}{.7pt}
\begin{picture}(60,40)(0,0)
\put(30,0){\line(0,1){10}}
\put(30,10){\line(-1,1){30}}
\put(30,10){\line(1,1){30}}
\put(15,25){\line(1,1){15}}
\end{picture}}%\kern 0.4ex}

\def\arbreC{%\kern-0.4ex
\setlength{\unitlength}{.7pt}
\begin{picture}(60,40)(0,0)
\put(30,0){\line(0,1){10}}
\put(30,10){\line(-1,1){30}}
\put(30,10){\line(1,1){30}}
\put(45,25){\line(-1,1){15}}
\end{picture}}%\kern 0.4ex}

\def\arbreBC{%\kern-0.4ex
\setlength{\unitlength}{.7pt}
\begin{picture}(60,40)(0,0)
\put(30,0){\line(0,1){40}}
\put(30,10){\line(-1,1){30}}
\put(30,10){\line(1,1){30}}
\end{picture}}%\kern 0.4ex}

\def\arbreBCd{%\kern-0.4ex
\setlength{\unitlength}{.7pt}
\begin{picture}(60,40)(0,0)
\put(30,0){\line(0,1){40}}
\put(30,10){\line(-1,1){30}}
\put(30,10){\line(1,1){30}}
\put(40,5){\small{(x,y)}}
\end{picture}}%\kern 0.4ex}

\def\arbreun{%\kern-0.4ex
\setlength{\unitlength}{.7pt}
\begin{picture}(60,40)(0,0)
\put(30,0){\line(0,1){10}}
\put(30,10){\line(-1,1){30}}
\put(30,10){\line(1,1){30}}
\put(20,20){\line(1,1){20}}
\put(10,30){\line(1,1){10}}
\end{picture}}%\kern 0.4ex}

\def\arbredeux{%\kern-0.4ex
\setlength{\unitlength}{.7pt}
\begin{picture}(60,40)(0,0)
\put(30,0){\line(0,1){10}}
\put(30,10){\line(-1,1){30}}
\put(30,10){\line(1,1){30}}
\put(20,20){\line(1,1){20}}
\put(30,30){\line(-1,1){10}}
\end{picture}}%\kern 0.4ex}

\def\arbretrois{%\kern-0.4ex
\setlength{\unitlength}{1pt}
\begin{picture}(60,40)(0,0)
\put(30,0){\line(0,1){10}}
\put(30,10){\line(-1,1){30}}
\put(30,10){\line(1,1){30}}
\put(50,30){\line(-1,1){10}}
\put(10,30){\line(1,1){10}}
\end{picture}}%\kern 0.4ex}

\def\arbretroisd{%\kern-0.4ex
\setlength{\unitlength}{1pt}
\begin{picture}(60,40)(0,0)
\put(30,0){\line(0,1){10}}
\put(30,10){\line(-1,1){30}}
\put(30,10){\line(1,1){30}}
\put(50,30){\line(-1,1){10}}
\put(10,30){\line(1,1){10}}
\put(2,27){\small{x}} %(x,y)
\put(33,5){\small{y}}
\put(55,27){\small{z}} %(x,y)
\end{picture}}%\kern 0.4ex}

\def\arbrequatre{%\kern-0.4ex
\setlength{\unitlength}{.7pt}
\begin{picture}(60,40)(0,0)
\put(30,0){\line(0,1){10}}
\put(30,10){\line(-1,1){30}}
\put(30,10){\line(1,1){30}}
\put(40,20){\line(-1,1){20}}
\put(30,30){\line(1,1){10}}
\end{picture}}%\kern 0.4ex}

\def\arbrecinq{%\kern-0.4ex
\setlength{\unitlength}{.7pt}
\begin{picture}(60,40)(0,0)
\put(30,0){\line(0,1){10}}
\put(30,10){\line(-1,1){30}}
\put(30,10){\line(1,1){30}}
\put(40,20){\line(-1,1){20}}
\put(50,30){\line(-1,1){10}}
\end{picture}}%\kern 0.4ex}

\def\arbreuut{%\kern-0.4ex
\setlength{\unitlength}{.7pt}
\begin{picture}(60,40)(0,0)
\put(30,0){\line(0,1){10}}
\put(30,10){\line(-1,1){30}}
\put(30,10){\line(1,1){30}}
\put(20,20){\line(0,1){20}}
\put(20,20){\line(1,1){20}}
\end{picture}}%\kern 0.4ex}

\def\arbretut{%\kern-0.4ex
\setlength{\unitlength}{.7pt}
\begin{picture}(60,40)(0,0)
\put(30,0){\line(0,1){10}}
\put(30,10){\line(-1,1){30}}
\put(30,10){\line(1,1){30}}
\put(30,10){\line(0,1){20}}
\put(30,30){\line(1,1){10}}
\put(30,30){\line(-1,1){10}}
\end{picture}}%\kern 0.4ex}

\def\arbretuu{%\kern-0.4ex
\setlength{\unitlength}{.7pt}
\begin{picture}(60,40)(0,0)
\put(30,0){\line(0,1){10}}
\put(30,10){\line(-1,1){30}}
\put(30,10){\line(1,1){30}}
\put(40,20){\line(0,1){20}}
\put(40,20){\line(-1,1){20}}
\end{picture}}%\kern 0.4ex}

\def\arbreutt{%\kern-0.4ex
\setlength{\unitlength}{.7pt}
\begin{picture}(60,40)(0,0)
\put(30,0){\line(0,1){10}}
\put(30,10){\line(-1,1){30}}
\put(30,10){\line(1,1){30}}
\put(10,30){\line(1,1){10}}
\put(30,10){\line(0,1){30}}
\end{picture}}%\kern 0.4ex}

\def\arbrettu{%\kern-0.4ex
\setlength{\unitlength}{.7pt}
\begin{picture}(60,40)(0,0)
\put(30,0){\line(0,1){10}}
\put(30,10){\line(-1,1){30}}
\put(30,10){\line(1,1){30}}
\put(50,30){\line(-1,1){10}}
\put(30,10){\line(0,1){30}}
\end{picture}}%\kern 0.4ex}

\def\arbrettt{%\kern-0.4ex
\setlength{\unitlength}{.7pt}
\begin{picture}(60,40)(0,0)
\put(30,0){\line(0,1){10}}
\put(30,10){\line(-1,1){30}}
\put(30,10){\line(1,1){30}}
\put(30,10){\line(-1,2){15}}
\put(30,10){\line(1,2){15}}
\end{picture}}%\kern 0.4ex}

%%%%%%%%small-scale trees
\def\s1tree{\!\!\includegraphics[scale=0.41]{1tree.eps}}
%%%%trees with 3 leaves
\def\sa2tree{\!\!\includegraphics[scale=0.41]{2tree.eps}}
\def\sb3tree{\!\!\includegraphics[scale=0.41]{3tree.eps}}
\def\sc4tree{\!\!\includegraphics[scale=0.41]{4tree.eps}}

%%%%trees with 2 leaves
\def\1tree{\!\!\includegraphics[scale=0.51]{1tree.eps}}
%%%%trees with 3 leaves
\def\2tree{\!\!\includegraphics[scale=0.51]{2tree.eps}}
\def\3tree{\!\!\includegraphics[scale=0.51]{3tree.eps}}
\def\4tree{\!\!\includegraphics[scale=0.51]{4tree.eps}}
%%%%trees with 4 leaves
\def\5tree{\!\!\includegraphics[scale=0.51]{5tree.eps}}
\def\6tree{\!\!\includegraphics[scale=0.51]{6tree.eps}}
\def\7tree{\!\!\includegraphics[scale=0.51]{7tree.eps}}
\def\8tree{\!\!\includegraphics[scale=0.51]{8tree.eps}}
\def\9tree{\!\!\includegraphics[scale=0.51]{9tree.eps}}
\def\a1tree{\!\!\includegraphics[scale=0.51]{10tree.eps}}
\def\b1tree{\!\!\includegraphics[scale=0.51]{11tree.eps}}
\def\c1tree{\!\!\includegraphics[scale=0.51]{12tree.eps}}
\def\d1tree{\!\!\includegraphics[scale=0.51]{13tree.eps}}
\def\e1tree{\!\!\includegraphics[scale=0.51]{14tree.eps}}
\def\f1tree{\!\!\includegraphics[scale=0.51]{15tree.eps}}
%%%%%decorated trees
\def\dec1tree{\!\!\includegraphics[scale=0.51]{dectree1.eps}}
\def\xtree{\!\!\includegraphics[scale=0.41]{xtree.eps}}
\def\xyztree{\!\!\includegraphics[scale=0.41]{xyztree.eps}}

%%%%%%%%%%%%%%%%%%% EGG trees

\def\ta1{{\scalebox{0.25}{ %%%%%%%%%%%%%%%%%%%%%%%%%%%%%%%%%\ta1
\begin{picture}(12,12)(38,-38)
\SetWidth{0.5} \SetColor{Black} \Vertex(45,-33){5.66}
\end{picture}}}}

\def\tb2{{\scalebox{0.25}{ %%%%%%%%%%%%%%%%%%%%%%%%%%%%%%%%%\tb2
\begin{picture}(12,42)(38,-38)
\SetWidth{0.5} \SetColor{Black} \Vertex(45,-3){5.66}
\SetWidth{1.0} \Line(45,-3)(45,-33) \SetWidth{0.5}
\Vertex(45,-33){5.66}
\end{picture}}}}

\def\tc3{{\scalebox{0.25}{ %%%%%%%%%%%%%%%%%%%%%%%%%%%%%%%%%\tc3
\begin{picture}(12,72)(38,-38)
\SetWidth{0.5} \SetColor{Black} \Vertex(45,27){5.66}
\SetWidth{1.0} \Line(45,27)(45,-3) \SetWidth{0.5}
\Vertex(45,-33){5.66} \SetWidth{1.0} \Line(45,-3)(45,-33)
\SetWidth{0.5} \Vertex(45,-3){5.66}
\end{picture}}}}

\def\td31{{\scalebox{0.25}{ %%%%%%%%%%%%%%%%%%%%%%%%%%%%%%%%%\td31
\begin{picture}(42,42)(23,-38)
\SetWidth{0.5} \SetColor{Black} \Vertex(45,-3){5.66}
\Vertex(30,-33){5.66} \Vertex(60,-33){5.66} \SetWidth{1.0}
\Line(45,-3)(30,-33) \Line(60,-33)(45,-3)
\end{picture}}}}

\def\te4{{\scalebox{0.25}{ %%%%%%%%%%%%%%%%%%%%%%%%%%%%%%%%%\te4
\begin{picture}(12,102)(38,-8)
\SetWidth{0.5} \SetColor{Black} \Vertex(45,57){5.66}
\Vertex(45,-3){5.66} \Vertex(45,27){5.66} \Vertex(45,87){5.66}
\SetWidth{1.0} \Line(45,57)(45,27) \Line(45,-3)(45,27)
\Line(45,57)(45,87)
\end{picture}}}}

\def\tf41{{\scalebox{0.25}{ %%%%%%%%%%%%%%%%%%%%%%%%%%%%%%%%%\tf41
\begin{picture}(42,72)(38,-8)
\SetWidth{0.5} \SetColor{Black} \Vertex(45,27){5.66}
\Vertex(45,-3){5.66} \SetWidth{1.0} \Line(45,27)(45,-3)
\SetWidth{0.5} \Vertex(60,57){5.66} \SetWidth{1.0}
\Line(45,27)(60,57) \SetWidth{0.5} \Vertex(75,27){5.66}
\SetWidth{1.0} \Line(75,27)(60,57)
\end{picture}}}}

\def\tg42{{\scalebox{0.25}{ %%%%%%%%%%%%%%%%%%%%%%%%%%%%%%%%%\tg42
\begin{picture}(42,72)(8,-8)
\SetWidth{0.5} \SetColor{Black} \Vertex(45,27){5.66}
\Vertex(45,-3){5.66} \SetWidth{1.0} \Line(45,27)(45,-3)
\SetWidth{0.5} \Vertex(15,27){5.66} \Vertex(30,57){5.66}
\SetWidth{1.0} \Line(15,27)(30,57) \Line(45,27)(30,57)
\end{picture}}}}

\def\th43{{\scalebox{0.25}{ %%%%%%%%%%%%%%%%%%%%%%%%%%%%%%%%%\th43
\begin{picture}(42,42)(8,-8)
\SetWidth{0.5} \SetColor{Black} \Vertex(45,-3){5.66}
\Vertex(15,-3){5.66} \Vertex(30,27){5.66} \SetWidth{1.0}
\Line(15,-3)(30,27) \Line(45,-3)(30,27) \Line(30,27)(30,-3)
\SetWidth{0.5} \Vertex(30,-3){5.66}
\end{picture}}}}

\def\thII43{{\scalebox{0.25}{ %%%%%%%%%%%%%%%%%%%%%%%%%%%%%%%%%\th43
\begin{picture}(72,57) (68,-128)
    \SetWidth{0.5}
    \SetColor{Black}
    \Vertex(105,-78){5.66}
    \SetWidth{1.5}
    \Line(105,-78)(75,-123)
    \Line(105,-78)(105,-123)
    \Line(105,-78)(135,-123)
    \SetWidth{0.5}
    \Vertex(75,-123){5.66}
    \Vertex(105,-123){5.66}
    \Vertex(135,-123){5.66}
  \end{picture}
  }}}

\def\thj44{{\scalebox{0.25}{ %%%%%%%%%%%%%%%%%%%%%%%%%%%%%%%%%\thj44
\begin{picture}(42,72)(8,-8)
\SetWidth{0.5} \SetColor{Black} \Vertex(30,57){5.66}
\SetWidth{1.0} \Line(30,57)(30,27) \SetWidth{0.5}
\Vertex(30,27){5.66} \SetWidth{1.0} \Line(45,-3)(30,27)
\SetWidth{0.5} \Vertex(45,-3){5.66} \Vertex(15,-3){5.66}
\SetWidth{1.0} \Line(15,-3)(30,27)
\end{picture}}}}

\def\ti5{{\scalebox{0.25}{ %%%%%%%%%%%%%%%%%%%%%%%%%%%%%%%%%\ti5
\begin{picture}(12,132)(23,-8)
\SetWidth{0.5} \SetColor{Black} \Vertex(30,117){5.66}
\SetWidth{1.0} \Line(30,117)(30,87) \SetWidth{0.5}
\Vertex(30,87){5.66} \Vertex(30,57){5.66} \Vertex(30,27){5.66}
\Vertex(30,-3){5.66} \SetWidth{1.0} \Line(30,-3)(30,27)
\Line(30,27)(30,57) \Line(30,87)(30,57)
\end{picture}}}}

\def\tj51{{\scalebox{0.25}{ %%%%%%%%%%%%%%%%%%%%%%%%%%%%%%%%%\tj51
\begin{picture}(42,102)(53,-38)
\SetWidth{0.5} \SetColor{Black} \Vertex(61,27){4.24}
\SetWidth{1.0} \Line(75,57)(90,27) \Line(60,27)(75,57)
\SetWidth{0.5} \Vertex(90,-3){5.66} \Vertex(60,27){5.66}
\Vertex(75,57){5.66} \Vertex(90,-33){5.66} \SetWidth{1.0}
\Line(90,-33)(90,-3) \Line(90,-3)(90,27) \SetWidth{0.5}
\Vertex(90,27){5.66}
\end{picture}}}}

\def\tk52{{\scalebox{0.25}{ %%%%%%%%%%%%%%%%%%%%%%%%%%%%%%%%%\tk52
\begin{picture}(42,102)(23,-8)
\SetWidth{0.5} \SetColor{Black} \Vertex(60,57){5.66}
\Vertex(45,87){5.66} \SetWidth{1.0} \Line(45,87)(60,57)
\SetWidth{0.5} \Vertex(30,57){5.66} \SetWidth{1.0}
\Line(30,57)(45,87) \SetWidth{0.5} \Vertex(30,-3){5.66}
\SetWidth{1.0} \Line(30,-3)(30,27) \SetWidth{0.5}
\Vertex(30,27){5.66} \SetWidth{1.0} \Line(30,57)(30,27)
\end{picture}}}}

\def\tl53{{\scalebox{0.25}{ %%%%%%%%%%%%%%%%%%%%%%%%%%%%%%%%%\tl53
\begin{picture}(42,102)(8,-8)
\SetWidth{0.5} \SetColor{Black} \Vertex(30,57){5.66}
\Vertex(30,27){5.66} \SetWidth{1.0} \Line(30,57)(30,27)
\SetWidth{0.5} \Vertex(30,87){5.66} \SetWidth{1.0}
\Line(30,27)(45,-3) \SetWidth{0.5} \Vertex(15,-3){5.66}
\SetWidth{1.0} \Line(15,-3)(30,27) \Line(30,57)(30,87)
\SetWidth{0.5} \Vertex(45,-3){5.66}
\end{picture}}}}

\def\tm54{{\scalebox{0.25}{ %%%%%%%%%%%%%%%%%%%%%%%%%%%%%%%%%\tm54
\begin{picture}(42,72)(8,-38)
\SetWidth{0.5} \SetColor{Black} \Vertex(30,-3){5.66}
\SetWidth{1.0} \Line(30,27)(30,-3) \Line(30,-3)(45,-33)
\SetWidth{0.5} \Vertex(15,-33){5.66} \SetWidth{1.0}
\Line(15,-33)(30,-3) \SetWidth{0.5} \Vertex(45,-33){5.66}
\SetWidth{1.0} \Line(30,-33)(30,-3) \SetWidth{0.5}
\Vertex(30,-33){5.66} \Vertex(30,27){5.66}
\end{picture}}}}

\def\tn55{{\scalebox{0.25}{ %%%%%%%%%%%%%%%%%%%%%%%%%%%%%%%%%\tn55
\begin{picture}(42,72)(8,-38)
\SetWidth{0.5} \SetColor{Black} \Vertex(15,-33){5.66}
\Vertex(45,-33){5.66} \Vertex(30,27){5.66} \SetWidth{1.0}
\Line(45,-33)(45,-3) \SetWidth{0.5} \Vertex(45,-3){5.66}
\Vertex(15,-3){5.66} \SetWidth{1.0} \Line(30,27)(45,-3)
\Line(15,-3)(30,27) \Line(15,-3)(15,-33)
\end{picture}}}}

\def\tp56{{\scalebox{0.25}{ %%%%%%%%%%%%%%%%%%%%%%%%%%%%%%%%%\tp56
\begin{picture}(66,111)(0,0)
\SetWidth{0.5} \SetColor{Black} \Vertex(30,66){5.66}
\Vertex(45,36){5.66} \SetWidth{1.0} \Line(30,66)(45,36)
\Line(15,36)(30,66) \SetWidth{0.5} \Vertex(30,6){5.66}
\Vertex(60,6){5.66} \SetWidth{1.0} \Line(60,6)(45,36)
\SetWidth{0.5}
\SetWidth{1.0} \Line(45,36)(30,6) \SetWidth{0.5}
\Vertex(15,36){5.66}
\end{picture}}}}

\def\tq57{{\scalebox{0.25}{ %%%%%%%%%%%%%%%%%%%%%%%%%%%%%%%%%\tq57
\begin{picture}(81,111)(0,0)
\SetWidth{0.5} \SetColor{Black} \Vertex(45,36){5.66}
\Vertex(30,6){5.66} \Vertex(60,6){5.66} \SetWidth{1.0}
\Line(60,6)(45,36) \SetWidth{0.5}
\SetWidth{1.0} \Line(45,36)(30,6) \SetWidth{0.5}
\Vertex(75,36){5.66} \SetWidth{1.0} \Line(45,36)(60,66)
\Line(60,66)(75,36) \SetWidth{0.5} \Vertex(60,66){5.66}
\end{picture}}}}

\def\tr58{{\scalebox{0.25}{ %%%%%%%%%%%%%%%%%%%%%%%%%%%%%%%%%\tr58
\begin{picture}(81,111)(0,0)
\SetWidth{0.5} \SetColor{Black} \Vertex(60,6){5.66}
\Vertex(75,36){5.66} \SetWidth{1.0} \Line(60,66)(75,36)
\SetWidth{0.5} \Vertex(60,66){5.66}
\SetWidth{1.0} \Line(60,36)(60,66) \Line(60,6)(60,36)
\SetWidth{0.5} \Vertex(60,36){5.66} \Vertex(45,36){5.66}
\SetWidth{1.0} \Line(60,66)(45,36)
\end{picture}}}}

\def\ts59{{\scalebox{0.25}{ %%%%%%%%%%%%%%%%%%%%%%%%%%%%%%%%%\ts59
\begin{picture}(81,111)(0,0)
\SetWidth{0.5} \SetColor{Black}
\Vertex(75,36){5.66} \SetWidth{1.0} \Line(60,66)(75,36)
\SetWidth{0.5} \Vertex(60,66){5.66}
\SetWidth{1.0} \Line(60,36)(60,66) \SetWidth{0.5}
\Vertex(60,36){5.66} \Vertex(45,36){5.66} \SetWidth{1.0}
\Line(60,66)(45,36) \Line(75,6)(75,36) \SetWidth{0.5}
\Vertex(75,6){5.66}
\end{picture}}}}

\def\tt591{{\scalebox{0.25}{ %%%%%%%%%%%%%%%%%%%%%%%%%%%%%%%%%\tt591
\begin{picture}(81,111)(0,0)
\SetWidth{0.5} \SetColor{Black}
\Vertex(75,36){5.66} \SetWidth{1.0} \Line(60,66)(75,36)
\SetWidth{0.5} \Vertex(60,66){5.66}
\SetWidth{1.0} \Line(60,36)(60,66) \SetWidth{0.5}
\Vertex(60,36){5.66} \Vertex(45,36){5.66} \SetWidth{1.0}
\Line(60,66)(45,36) \SetWidth{0.5} \Vertex(45,6){5.66}
\SetWidth{1.0} \Line(45,6)(45,36)
\end{picture}}}}

\def\treeUNO{{{\scalebox{0.25}{ %%%%%%%%%%%%%%%%%%%%%%%%%%%\treeUNO
 \begin{picture}(72,72) (83,-98)
    \SetWidth{0.5}
    \SetColor{Black}
    \Vertex(105,-63){5.66}
    \Vertex(90,-93){5.66}
    \Vertex(120,-93){5.66}
    \SetWidth{1.1}
    \Line(105,-63)(90,-93)
    \Line(105,-63)(120,-93)
    \SetWidth{0.5}
    \Vertex(150,-93){5.66}
    \Vertex(120,-33){5.66}
    \SetWidth{1.5}
    \Line(120,-33)(150,-93)
    \Line(120,-33)(105,-63)
  \end{picture}}}}}

\def\treeDUO{{{\scalebox{0.25}{ %%%%%%%%%%%%%%%%%%%%%%%%%%%\treeDUE
 \begin{picture}(102,72) (53,-98)
    \SetWidth{0.5}
    \SetColor{Black}
    \Vertex(105,-63){5.66}
    \Vertex(90,-93){5.66}
    \Vertex(120,-93){5.66}
    \SetWidth{1.5}
    \Line(105,-63)(90,-93)
    \Line(105,-63)(120,-93)
    \SetWidth{0.5}
    \Vertex(105,-33){5.66}
    \SetWidth{1.5}
    \Line(105,-33)(105,-63)
    \SetWidth{0.5}
    \Vertex(150,-93){5.66}
    \Vertex(60,-93){5.66}
    \SetWidth{1.5}
    \Line(105,-33)(60,-93)
    \Line(105,-33)(150,-93)
  \end{picture}}}}}

%%%%%%%%%%%%%%%%%%%%%%%%%%%%%%%%%%%%%%%%%%%%%%%%%%%%%%%%%%%%%%%Decorated trees%

\def\ydec31{\!\!\includegraphics[scale=0.5]{ydec31.eps}}

\def\kyldec31{\!\!\includegraphics[scale=0.5]{kyldec31.eps}}

\def\yldec31{\!\!\includegraphics[scale=0.5]{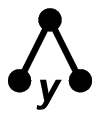}}

\def\xldec41r{\!\!\includegraphics[scale=0.5]{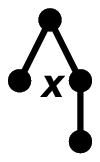}}

\def\xyldec43{\!\!\includegraphics[scale=0.5]{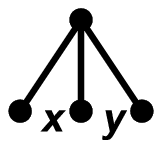}}

\def\xtd31{\!\!\includegraphics[scale=1]{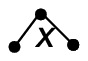}}

\def\xthj44{\!\!\includegraphics[scale=1]{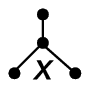}}

%\def\treeUNO{\includegraphics[scale=0.4]{treeUNO.eps}}
%\def\treeDUO{\includegraphics[scale=0.4]{treeDUO.eps}}
%\def\treeTRES{\includegraphics[scale=0.5]{treeTRES.eps}}

%%%%%%%%%%%%%%%%%%%%%%%%%%%%%%%%%%%%%%%%%%%%%%%%%%%%%%%%%%%%%%%%%%%%%%%%%%%%

\def\lta1{{\scalebox{0.25}{ %%%%%%%%%%%%%%%%%%%%%%%%%%%%%%%%%\lta1
\begin{picture}(0,45)(60,-15)
\SetWidth{1.5} \SetColor{Black} \Line(60,30)(60,-15)
\end{picture}}}}

\def\ltb2{{\scalebox{0.25}{ %%%%%%%%%%%%%%%%%%%%%%%%%%%%%%%%%\ltb2
\begin{picture}(12,45)(53,-15)
\SetWidth{1.5} \SetColor{Black} \Line(60,30)(60,-15)
\SetWidth{0.5} \Vertex(60,0){5.66}
\end{picture}}}}

\def\ltc3{{\scalebox{0.25}{ %%%%%%%%%%%%%%%%%%%%%%%%%%%%%%%%%\ltc3
\begin{picture}(12,75)(53,-15)
\SetWidth{0.5} \SetColor{Black} \Vertex(60,30){5.66}
\SetWidth{1.5} \Line(60,60)(60,-15) \SetWidth{0.5}
\Vertex(60,0){5.66}
\end{picture}}}}

\def\ltd31{{\scalebox{0.25}{ %%%%%%%%%%%%%%%%%%%%%%%%%%%%%%%%%\ltd31
\begin{picture}(75,90)(0,0)
\SetWidth{0.5} \SetColor{Black}
\Vertex(60,15){5.66} \SetWidth{1.5} \Line(45,45)(60,15)
\Line(60,15)(75,45) \Line(60,15)(60,0)
\end{picture}}}}

\delete{
\begin{picture}(75,90)(0,0)
\SetWidth{0.5} \SetColor{Black} \% \Vertex(60,15){5.66}
\SetWidth{1.5} \Line(45,45)(60,15) \Line(60,15)(75,45)
\Line(60,15)(60,0)
\end{picture}}

\def\lte4{{\scalebox{0.25}{ %%%%%%%%%%%%%%%%%%%%%%%%%%%%%%%%%\lte4
\begin{picture}(66,120)(0,0)
\SetWidth{0.5} \SetColor{Black}
\Vertex(60,45){5.66} \Vertex(60,75){5.66} \SetWidth{1.5}
\Line(60,105)(60,0) \SetWidth{0.5} \Vertex(60,15){5.66}
\end{picture}}}}

\def\ltf41l{{\scalebox{0.25}{ %%%%%%%%%%%%%%%%%%%%%%%%%%%%%%%%%\ltf41l
\begin{picture}(75,120)(0,0)
\SetWidth{0.5} \SetColor{Black}
\Vertex(60,15){5.66} \SetWidth{1.5} \Line(60,0)(60,15)
\Line(60,15)(45,45) \Line(60,15)(75,45) \SetWidth{0.5}
\Vertex(45,45){5.66} \SetWidth{1.5} \Line(45,45)(45,75)
\end{picture}}}}

\def\ltg41r{{\scalebox{0.25}{ %%%%%%%%%%%%%%%%%%%%%%%%%%%%%%%%%\ltg41r
\begin{picture}(81,120)(0,0)
\SetWidth{0.5} \SetColor{Black}
\Vertex(60,15){5.66} \SetWidth{1.5} \Line(75,45)(75,75)
\Line(60,0)(60,15) \Line(60,15)(45,45) \Line(60,15)(75,45)
\SetWidth{0.5} \Vertex(75,45){5.66}
\end{picture}}}}

\delete{
\begin{picture}(81,120)(0,0)
\SetWidth{0.5} \SetColor{Black}
\Vertex(60,15){5.66} \SetWidth{1.5} \Line(75,45)(75,75)
\Line(60,0)(60,15) \Line(60,15)(45,45) \Line(60,15)(75,45)
\SetWidth{0.5} \Vertex(75,45){5.66}
\end{picture}}

\def\lth42{{\scalebox{0.25}{ %%%%%%%%%%%%%%%%%%%%%%%%%%%%%%%%%\lth42
\begin{picture}(75,150)(0,0)
\SetWidth{0.5} \SetColor{Black}
\Vertex(60,45){5.66} \SetWidth{1.5} \Line(60,30)(60,45)
\Line(60,45)(45,75) \Line(60,45)(75,75) \SetWidth{0.5}
\Vertex(60,15){5.66} \SetWidth{1.5} \Line(60,0)(60,30)
\end{picture}}}}

\def\lti43{{\scalebox{0.25}{ %%%%%%%%%%%%%%%%%%%%%%%%%%%%%%%%%\lti43
\begin{picture}(75,150)(0,0)
\SetWidth{0.5} \SetColor{Black}
\SetWidth{1.5} \Line(60,30)(60,45) \SetWidth{0.5}
\Vertex(60,15){5.66} \SetWidth{1.5} \Line(60,0)(60,30)
\Line(60,15)(75,45) \Line(60,15)(45,45)
\end{picture}}}}

\delete{
\begin{picture}(75,150)(0,0)
\SetWidth{0.5} \SetColor{Black}
\SetWidth{1.5} \Line(60,30)(60,45) \SetWidth{0.5}
\Vertex(60,15){5.66} \SetWidth{1.5} \Line(60,0)(60,30)
\Line(60,15)(75,45) \Line(60,15)(45,45)
\end{picture}}

\delete{
%%%%%%%%%%%%%%%%%%%%%%%%%%%%%%%%%%%%%%%%%%%%%%%%%%%%%%%%%%%%%%%%%%%%%%%%%%%%

\def\lta1{{\scalebox{0.25}{ %%%%%%%%%%%%%%%%%%%%%%%%%%%%%%%%%\lta1
\begin{picture}(0,45)(60,-15)
\SetWidth{1.5} \SetColor{Black} \Line(60,30)(60,-15)
\end{picture}}}}

\def\ltb2{{\scalebox{0.25}{ %%%%%%%%%%%%%%%%%%%%%%%%%%%%%%%%%\ltb2
\begin{picture}(12,45)(53,-15)
\SetWidth{1.5} \SetColor{Black} \Line(60,30)(60,-15)
\SetWidth{0.5} \Vertex(60,0){5.66}
\end{picture}}}}

\def\ltc3{{\scalebox{0.25}{ %%%%%%%%%%%%%%%%%%%%%%%%%%%%%%%%%\ltc3
\begin{picture}(12,75)(53,-15)
\SetWidth{0.5} \SetColor{Black} \Vertex(60,30){5.66}
\SetWidth{1.5} \Line(60,60)(60,-15) \SetWidth{0.5}
\Vertex(60,0){5.66}
\end{picture}}}}

\def\ltd31{{\scalebox{0.25}{ %%%%%%%%%%%%%%%%%%%%%%%%%%%%%%%%%\ltd31
\begin{picture}(75,90)(0,0)
\SetWidth{0.5} \SetColor{Black}
\Vertex(60,15){5.66} \SetWidth{1.5} \Line(45,45)(60,15)
\Line(60,15)(75,45) \Line(60,15)(60,0)
\end{picture}}}}

\def\lte4{{\scalebox{0.25}{ %%%%%%%%%%%%%%%%%%%%%%%%%%%%%%%%%\lte4
\begin{picture}(66,120)(0,0)
\SetWidth{0.5} \SetColor{Black}
\Vertex(60,45){5.66} \Vertex(60,75){5.66} \SetWidth{1.5}
\Line(60,105)(60,0) \SetWidth{0.5} \Vertex(60,15){5.66}
\end{picture}}}}

\def\ltf41l{{\scalebox{0.25}{ %%%%%%%%%%%%%%%%%%%%%%%%%%%%%%%%%\ltf41l
\begin{picture}(75,120)(0,0)
\SetWidth{0.5} \SetColor{Black}
\Vertex(60,15){5.66} \SetWidth{1.5} \Line(60,0)(60,15)
\Line(60,15)(45,45) \Line(60,15)(75,45) \SetWidth{0.5}
\Vertex(45,45){5.66} \SetWidth{1.5} \Line(45,45)(45,75)
\end{picture}}}}

\def\ltg41r{{\scalebox{0.25}{ %%%%%%%%%%%%%%%%%%%%%%%%%%%%%%%%%\ltg41r
\begin{picture}(81,120)(0,0)
\SetWidth{0.5} \SetColor{Black}
\Vertex(60,15){5.66} \SetWidth{1.5} \Line(75,45)(75,75)
\Line(60,0)(60,15) \Line(60,15)(45,45) \Line(60,15)(75,45)
\SetWidth{0.5} \Vertex(75,45){5.66}
\end{picture}}}}

\def\lth42{{\scalebox{0.25}{ %%%%%%%%%%%%%%%%%%%%%%%%%%%%%%%%%\lth42
\begin{picture}(75,150)(0,0)
\SetWidth{0.5} \SetColor{Black}
\Vertex(60,45){5.66} \SetWidth{1.5} \Line(60,30)(60,45)
\Line(60,45)(45,75) \Line(60,45)(75,75) \SetWidth{0.5}
\Vertex(60,15){5.66} \SetWidth{1.5} \Line(60,0)(60,30)
\end{picture}}}}

\def\lti43{{\scalebox{0.25}{ %%%%%%%%%%%%%%%%%%%%%%%%%%%%%%%%%\lti43
\begin{picture}(75,150)(0,0)
\SetWidth{0.5} \SetColor{Black}
\SetWidth{1.5} \Line(60,30)(60,45) \SetWidth{0.5}
\Vertex(60,15){5.66} \SetWidth{1.5} \Line(60,0)(60,30)
\Line(60,15)(75,45) \Line(60,15)(45,45)
\end{picture}}}}
}
%%%%%%%%%%%%%%%%%%%%%%%%%%%%%%%%%%%%%%%%%%%%%%%%%%%%%%%%%%%%%%%%%%%%%%%%%%%%

%%%%%%%%%%%%%%%%%%%%% leave trees
%%%%%%%%%%%%%%%%%%%%%%%%%%%%%%%%%%%%%%%%%%%%%%%%%%%%%%%%%%%

\def\lxtree{\includegraphics[scale=.7]{lxtree}}
\def\xdtree{\includegraphics[scale=.7]{xdtree}}
\def\xdtreend{\includegraphics[scale=.2]{xdtreend}} % xdtree no decoration
\def\xytree{\includegraphics[scale=.7]{xytree}}
\def\xytreend{\includegraphics[scale=.2]{xytreend}} % xytree no decoration

\def\xdtreed{\includegraphics[scale=.3]{xdtreed}}
\def\xytreed{\includegraphics[scale=.3]{xytreed}}
\def\lxdtree{\includegraphics[scale=.7]{lxdtree}}
\def\dtree{\includegraphics[scale=0.5]{dtree}}
\def\xyleft{\includegraphics[scale=.7]{xyleft}}
\def\xyright{\includegraphics[scale=.7]{xyright}}
\def\yzdtree{\includegraphics[scale=.7]{yzdtree}}
\def\xyztreea{\includegraphics[scale=.7]{xyztree1}}
\def\xyztreeb{\includegraphics[scale=.7]{xyztree2}}
\def\xyztreec{\includegraphics[scale=.7]{xyztree3}}
\def\tableps{\includegraphics[scale=1.4]{tableps}}

\def\ccxtree{\includegraphics[scale=1]{ccxtree}}
\def\ccxtreedec{\includegraphics[scale=1]{ccxtreedec}}
\def\ccetree{\includegraphics[scale=1]{ccetree}}
\def\ccllcrcr{\includegraphics[scale=1]{ccllcrcr}}
\def\cclclcrcr{\includegraphics[scale=.5]{cclclcrcr}}
\def\ccllxryr{\includegraphics[scale=.5]{ccllxryr}}
\def\cclxlyrzr{\includegraphics[scale=1]{cclxlyrzr}}
\def\ccIII{\includegraphics[scale=1]{ccIII}}

%%%%%%%%%%%%%%%%%%%%%%%%%%%%%%%%%%%%%%%%%%%%%%%%%%%%%%%%%%%%%%%%%%

\title{Free Rota--Baxter algebras and rooted trees}
\author{Kurusch Ebrahimi-Fard}
\address{I.H.\'E.S.
         Le Bois-Marie,
         35, Route de Chartres,
         F-91440 Bures-sur-Yvette,
         France}
         \email{kurusch@ihes.fr}
\author{Li Guo}
%\thanks{Corresponding author: Li Guo. Phone: 973-353-5156 ext. 30; Fax: 973-353-5270; E-mail: liguo@newark.rutgers.edu}
\address{Department of Mathematics and Computer Science,
         Rutgers University,
         Newark, NJ 07102}
\email{liguo@newark.rutgers.edu}

%\date{\today}

%\begin{document}

\begin{abstract}
A Rota--Baxter algebra, also known as a Baxter algebra, is an algebra with a linear operator satisfying a relation, called the Rota--Baxter relation, that generalizes the integration by parts formula. Most of the studies on Rota--Baxter algebras have been for commutative algebras. Two constructions of free commutative Rota--Baxter algebras were obtained by Rota and Cartier in the 1970s and a third one by Keigher and one of the authors in the 1990s in terms of mixable shuffles. Recently, noncommutative Rota--Baxter algebras have appeared both in physics in connection with the work of Connes and Kreimer on renormalization in perturbative quantum field theory, and in mathematics related to the work of Loday and Ronco on dendriform dialgebras and trialgebras.

This paper uses rooted trees and forests to give explicit constructions of free noncommutative Rota--Baxter algebras on modules and sets. This highlights the combinatorial nature of Rota--Baxter algebras and facilitates their further study. As an application, we obtain the unitarization of Rota--Baxter algebras.
\end{abstract}

% MSC: 16W99, 05C05, 08A50
% Keywords: Rota--Baxter algebras, free objects, rooted trees, planar rooted trees

\maketitle

%\tableofcontents

\setcounter{section}{0}
{\ }
\vspace{-1cm}

\section{Introduction}
\label{intro}
We construct the free Rota--Baxter algebra on a set $X$ in terms of angularly decorated rooted trees with $X$ as the decoration set. We also consider the more general case of free objects on modules. As an application, we prove the existence and uniqueness of the unitarization of Rota--Baxter algebras.
\medskip

%\subsection{Rota--Baxter algebras}

A Rota--Baxter algebra (also known as a Baxter algebra) is
an associative algebra $R$
with a linear endomorphism $P$ satisfying the {\bf{Rota--Baxter
relation}}:
\begin{equation}
    P(x)P(y) = P\big(P(x)y + xP(y) + \lambda xy\big),\ \forall\ x,y \in R.
    \mlabel{eq:RB}
\end{equation}
Here $\lambda$ is a fixed element in the base ring and is
sometimes denoted by $-\theta$. The relation was introduced by the
mathematician Glen E.~Baxter~\cite{Ba} in his probability study,
and was popularized mainly by the work of G.-C. Rota~\cite{Ro1, Ro2, R-S} and his school.

Note that the Rota--Baxter relation~(\mref{eq:RB}) is defined even if the binary operation is not associative. In fact,
such a relation for Lie algebras was
introduced independently by Belavin and
Drinfeld~\cite{B-D}, and Semenov-Tian-Shansky~\cite{STS1} in the
1980s, under the disguise of
$r$-matrices, of the (modified) classical Yang--Baxter equation,
named after the physicists Chen-ning Yang and Rodney Baxter.
Recently, there have been several interesting developments
of Rota--Baxter algebras in
theoretical physics and mathematics, including quantum field
theory~\cite{C-K1,C-K2,Kr1,Kr2,M-P1}, associative Yang--Baxter
equations~\cite{Ag1,Ag3}, shuffle
products~\cite{E-G1,G-K1,G-K2},
operads~\cite{A-L,EF1,E-G2,Le1,Le3}, Hopf algebras~\cite{A-G-K-O,E-G1}, combinatorics~\cite{Gu2} and number
theory~\cite{E-G6,Gu5,G-Z,M-P2}. The most prominent of these is the
work of Connes and Kreimer in their Hopf algebraic approach to
renormalization theory in perturbative quantum field
theory~\cite{C-K1,C-K2}, continued
in~\cite{E-G-G-V,E-G5,E-G-K2,E-G-K3}.
\medskip

%\subsection{Free commutative Rota--Baxter algebras}

Our goal in this paper is to give an explicit construction of
free {\em noncommutative} Rota--Baxter algebras in terms of
rooted trees. To help put this study in perspective, we briefly
review the interesting development of the {\em commutative} case.
Cartier~\cite{Ca} pointed out over thirty years ago ``The
existence of free (Rota--)Baxter algebras follows from well-known
arguments in universal algebra but remains quite immaterial as
long as the corresponding word problem is not solved in an
explicit way as Rota was the first to do." Both Rota's
aforementioned construction~\cite{Ro1} and the construction of
Cartier himself~\mcite{Ca} dealt with free
commutative Rota--Baxter algebras. Later, a third construction was
obtained by the second named author and Keigher~\cite{G-K1,G-K2}
as a generalization of shuffle product algebras.

These constructions of
free commutative Rota--Baxter algebras have
important implications. For example, Rota~\cite{Ro2,R-S} applied
his construction to give a proof of the celebrated Spitzer identity~\cite{Sp,E-G-K3} by relating it to
Waring's identity, another basic formula in combinatorics.
The product in Cartier's paper~\cite{Ca} is readily
seen to be the same as the one by Ehrenborg~\cite{Eh}
for  monomial quasi-symmetric functions and more
recently by Bradley~\cite{Br} to explicitly describe
stuffles and $q$-stuffles for multiple zeta values. Furthermore,
the mixable
shuffle product in the construction of~\cite{G-K1} appeared
also in the work of Goncharov~\cite{Go2} to study motivic
shuffle relations and the work of Hazewinkle~\cite{Ha} on overlapping shuffles. In~\cite{E-G1}, the mixable shuffle product
is shown to be the same as Hoffman's quasi-shuffle product~\cite{Ho}
which has played a fundamental role in the study of algebraic
relations among multiple zeta values.
There is also a description~\cite{A-H,Fa,Lo2} of quasi-shuffles in terms of piecewise linear paths (Delannoy paths).
\smallskip

Our consideration of the noncommutative case has motivations beyond a simple pursuit of generalization.
In the algebraic framework of Connes and Kreimer~\cite{C-K1,C-K2} for renormalization in quantum field theory, a regularized Feynman rule is viewed and studied as an algebra homomorphism from their Hopf algebra of Feynman diagrams to a Rota--Baxter algebra associated to the renormalization scheme. The renormalization and counter term for the Feynman rule are derived from the algebraic Birkhoff decomposition. In~\cite{E-G-K3}, the algebraic Birkhoff decomposition and the renormalization are shown to follow from the Atkinson decomposition and the Spitzer's identity in a noncommutative Rota--Baxter algebra. Since the Rota-Baxter algebra varies with the choice of a quantum field theory and renormalization scheme, it is desirable to investigate universal or free Rota-Baxter algebras.

In a more theoretical context, there have been quite strong
interests lately in possible noncommutative generalizations
of shuffles and quasi-shuffles (that is, mixable
shuffles). From the connection of these shuffles with
free commutative Rota--Baxter algebras mentioned above, such noncommutative
generalizations should be related to free noncommutative Rota--Baxter
algebras. Indeed one such generalization is the Hopf algebra of planar rooted trees of Loday and Ronco~\cite{L-R1} and we have shown in~\cite{E-G4} that
this algebra canonically embeds into a free noncommutative Rota--Baxter
algebra. The tree construction of free Rota--Baxter algebras obtained in this paper make this embedding a even more natural
tree-to-tree embedding. In fact, such an embedding has been our motivation to achieve a tree interpretation of free Rota--Baxter algebras.

It is also our hope that our explicit constructions of the free Rota--Baxter algebras here will lead to further studies of Rota--Baxter algebras. Indeed, some of such studies~\cite{A-M,G-S} have already been carried out concurrently with the writing of this paper. To compare with these and other related papers~\cite{E-G4,E-G0}, we note that there are different types of free Rota--Baxter algebras obtained from the adjoint functors of
the forgetful functors from the category of unitary Rota--Baxter algebras
to the categories of sets, modules, and algebras.
They give rise to free Rota--Baxter algebras generated by
(or on) a set, a module or an algebra. Further, by replacing unitary
algebras by nonunitary algebras, we get more forgetful
functors and their adjoint functors.
We summarize these categories and forgetful functors in the following diagram.
$$\xymatrix{ \mbox{Unitary Rota-Baxter algebras} \ar[r] \ar[d] &
\mbox{Unitary Algebras} \ar[r] \ar[d] &
\mbox{Modules} \ar[r] \ar^{=}[d] &
\mbox{Sets} \ar^{=}[d] \\
\mbox{Nonunitary Rota-Baxter algebras} \ar[r]  &
\mbox{Nonunitary Algebras} \ar[r]  &
\mbox{Modules} \ar[r] &
\mbox{Sets}
}
$$
The distinction between unitarity and nonunitarity for a Rota--Baxter algebra is more significant then for an associative algebra, because of the involvement of the Rota--Baxter
operator. In fact, it is with the help of
our constructions of unitary and nonunitary free Rota--Baxter
algebras that we prove the existence of unitarization of
Rota--Baxter algebras.
%\smallskip

In \cite{E-G4} the first construction of free Rota--Baxter algebras on another algebra were obtained in terms of bracketed words (called Rota--Baxter words). In the present paper we consider free Rota--Baxter algebras on a module and on a set in terms of rooted trees and forests. %
In~\cite{A-M}, free Rota--Baxter algebras are also constructed in terms of decorated rooted trees. It considers the singleton generating set
while we consider any generating set. The forms of the types of trees and decorations in the two papers are different with~\cite{A-M} using rooted tree with numerical decorations on the vertices and angles while us using rooted forests with angles decorated by the generating set or module. Also in~\cite{A-M} the Rota--Baxter algebras are constructed on the decorated trees while in our paper Rota--Baxter algebras are defined on forests without decoration and then are extended to forests with decorations.
Another related paper is~\cite{G-S} where enumeration, generating functions and algorithms of bracketed words in free Rota--Baxter algebras were studied. These aspects, in terms of trees and other combinatorial objects, were also considered in~\cite{A-M,Guop}.
\medskip

%\subsection{Outline}

This paper can be summarized by the following diagram of Rota--Baxter algebras.
$$
\xymatrix{
\bfk\, \calfo \ar@{_{(}->}[dd] & & & & \ncshao(M) \ar@{_{(}->}[dd]
&&& \ncshao(X) \ar@{_{(}->}[dd] & \ar^{\mbox{unitarization}}@{=>}[dd] \\
& \!\!\!\ar^{\mbox{angular} }_{\mbox{decoration}}@{=>}[rr] &&&
&\ar^{M=\bfk\, X}@{=>}[r] &&& \\
\bfk\, \calf  & & & & \ncsha(M) &&& \ncsha(X) &
}
$$
In Section~\mref{sec:rbt} we will consider the set of planar rooted forests $\calf$ and its subset $\calfo$ of ladder-free forests, and the corresponding free
$\bfk$-modules $\bfk\,\calf$ and $ \bfk\,\calfo$ over a commutative unitary ring $\bfk$. We equip these two modules with a Rota--Baxter algebra structure (Theorem~\mref{thm:freet} and Proposition~\mref{pp:freet0}).
By decorating angles of the forests in these Rota--Baxter algebras by elements of a module $M$, we construct in Section~\mref{sec:freem} the free unitary (resp. nonunitary) Rota--Baxter algebra
$\ncsha(M)$ (resp. $\ncshao(M)$) on $M$ in Theorem~\mref{thm:freem} (resp. Theorem~\mref{thm:freemo}).
By taking $M=\bfk\, X$ for a set $X$, we obtain free Rota--Baxter algebras on a set $X$ in Section~\mref{sec:freex} and display a canonical basis in the form of angularly decorated forests (Theorem~\mref{thm:freex}).
As an application of these free Rota--Baxter algebras,
the unitarization of Rota--Baxter algebras is obtained in
Section~\mref{sec:uni}.

\medskip

\noindent
{\bf Notations:}
In this paper, $\bfk$ is a commutative unitary ring. By a
$\bfk$-algebra we mean a unitary algebra over the base ring $\bfk$
unless otherwise stated. The same applies to Rota--Baxter
algebras. For a set $X$, let $\bfk\,X$ be the free $\bfk$-module
$\oplus_{x\in X} \bfk\,x$ generated by $X$. If $X$ is a semigroup
(resp. monoid), $\bfk\, X$ is equipped with the natural nonunitary (resp.
unitary) $\bfk$-algebra structure.
\medskip
%%%
%%%%%%%%%%%%%%%%%%%%%%%%%%%%%%%%%%%%%%%%%%%%%%%%%%%%%%%%%%%
%%%

\noindent{\bf Acknowledgements:}
The first named author
thanks the European Post-Doctoral Institute for a grant supporting
his stay at I.H.\'E.S. The second named
author is supported in part by NSF grant DMS 0505643 and a
Research Council grant from the Rutgers University. He thanks P. Cartier for helpful discussions. He also
thanks CIRM at Luminy where this work was started and thanks
Max-Planck Institute of Mathematics at Bonn where this work was completed.

%%%
%%%%%%%%%%%%%%%%%%%%%%%%%%%%%%%%%%%%%%%%%%%%%%%%%%%%%%%%%%%
%%%

\section{The Rota--Baxter algebra of planar rooted forests}
\mlabel{sec:rbt} We first obtain a Rota--Baxter algebra
structure on planar rooted forests and their various subsets. This allows us to give a uniform construction of free Rota--Baxter algebras in different settings in \S~\mref{sec:freem}. For other variations of this construction, see~\mcite{A-M,E-G0,E-G5}.
\subsection{Planar rooted forests}
For the convenience of the reader and for fixing notations, we
recall basic concepts and facts of planar rooted trees.
For references, see \cite{Di,We}.

A free tree is an undirected graph that is connected and
contains no cycles. A {\bf rooted tree} is a free tree in which a particular vertex has been distinguished as the {\bf root}. Such a distinguished vertex endows the tree with a directed graph structure when the edges of the tree are given the orientation of pointing away from the root. If two vertices of a rooted tree are connected by such an oriented edge, then the vertex on the side of the root is called the {\bf parent} and the vertex on the opposite side of the root is called a {\bf child}. A vertex with no children is called a {\bf leaf}. By our convention, in a tree with only one vertex, this vertex is a leaf, as well as the root. The number of edges in a path connecting two vertices in a rooted tree is called the {\bf length} of the path. The {\bf depth} $\depth(T)$ (or {\bf height}) of a rooted tree $T$ is the length of the longest path from its root to its leafs. A {\bf planar rooted tree} is a rooted tree with a fixed embedding into the plane.

There are two ways to
draw planar rooted trees. In one drawing all vertices are
represented by a dot and the root is usually at the top of the
tree. The following list shows the first few of them.
$$
 \ta1 \;\quad
 \tb2   \quad\;
 \tc3   \quad\;
 \td31  \quad\;
 \te4   \quad\;
 \tf41  \quad\;
 \tg42 \;\quad
 \thj44 \quad\;
 \th43 \;\quad
\delete{
  \ti5  \quad\;
 \tj51  \quad\;
 \tk52  \quad\;
 \tl53  \quad\;
 \tm54  \quad\;
 \tn55
% \ti5 \quad\; \tj51 \quad\; \ptj51 \quad\;
% \tk52 \quad\; \ptk52 \;\;\; \cdots \;\;\; \tp56 \;\quad \ptp56 \;\;\;
} \cdots
$$
Note that we distinguish the sides of the trees, so the
trees are planar. The tree $\onetree$ with only
the root is called the {\bf empty tree}.
This drawing is used, for example, in the above reference~\cite{Di,We} of trees and in the Hopf algebra of non-planar rooted trees of Connes and Kreimer~\cite{C-K0,C-K1}.

In the second drawing the leaf vertices are
removed with only the edges leading to them left, and the root,
placed at the bottom in opposite to the first drawing,  gets an
extra edge pointing down. The following list
shows the first few of them.
$$
 \lta1  \;\quad
 \ltb2  \;\quad
 \ltc3  \!\!\!\!\quad
 \ltd31 \!\!\!\!\quad
 \lte4   \!\!\!\!\quad
 \ltf41l \!\!\!\!\quad
 \ltg41r \!\!\!\!\quad
 \lth42 \!\!\!\!\quad
 \lti43 \quad \cdots
$$
This is used, for example in the Hopf algebra of planar rooted
trees of Loday and Ronco~\cite{Lo1,L-R1} and noncommutative
variation of the Connes-Kreimer Hopf algebra~\cite{Holt,Fo}.
In the following we will mostly use the first drawing.
%This usage is in contrast to the setting of
%Connes--Kreimer, where the empty tree is the empty set.
\medskip

%\subsubsection{Recursive definition of planar rooted trees}
Let $\calt$ be the set of planar rooted trees and let $\calf$ be
the free semigroup generated by $\calt$ in which the product is
denoted by $\sqcup$, called the concatenation. Thus each element in $\calf$ is a
noncommutative product $T_1\sqcup \cdots \sqcup T_n$ consisting of trees $T_1,\cdots, T_n\in\calt$, called a {\bf
planar rooted forest}.
We also use the abbreviation
\begin{equation}
T^{\sqcup n}=\underbrace{T\sqcup \cdots \sqcup T}_{n\ {\rm terms}}.
\mlabel{eq:sqpower}
\end{equation}
\begin{remark}
For the rest of this paper, a tree or forest means a planar rooted
one unless otherwise specified.
\end{remark}
We use the (grafting) {\bf brackets} $\lc T_1\sqcup \cdots \sqcup T_n\rc$ to denote the tree obtained by {\bf grafting}, that is,
by adding a new root together with an edge from the new root to the root of each of the trees $T_1,\cdots, T_n$. This is the $B^+$
operator in the work of Connes and Kreimer~\cite{C-K1}. The
operation is also denoted by $T_1 \vee \cdots \vee T_n$ in some
other literatures, such as in Loday and Ronco~\cite{Lo1,L-R1}.
Note that our
operation $\sqcup$ is different from $\vee$. Their relation is
$$
    \lc T_1 \sqcup \cdots \sqcup T_n\rc = T_1\vee \cdots \vee T_n.
$$
See~\mcite{Guop} for a general framework to view such algebraic structures with operators.

The {\bf depth} of a forest $F$ is the
maximal depth $\depth=\depth(F)$ of trees in $F$.
Clearly, $\depth(\lc F\rc)=\depth(F)+1$. The trees in a
forest $F$ are called root branches of $\lc F\rc$.
Furthermore, for a forest $F=T_1\sqcup \cdots \sqcup T_b$ with
trees $T_1,\cdots,T_b$, we define $b=b(F)$ to be the {\bf breadth}
of $F$. Let $\leaf(F)$ be the number of leafs of $F$. Then \begin{equation}
\leaf(F)=\sum_{i=1}^{b}\leaf(T_i).
\mlabel{eq:leaf}
\end{equation}

We will often use the following
recursive structure on forests.
For any subset $X$ of $\calf$, let $\sqmon{X}$ be the sub-semigroup
of $\calf$ generated by $X$.
Let $\calf_0=\sqmon{\onetree}$, consisting of forests
$\onetree^{\sqcup n}, n\geq 0$.
These are also the forests of depth zero.
Then recursively define
\begin{equation}
 \calf_n= \sqmon{\{\onetree\}\cup \lc \calf_{n-1}\rc}.
 \mlabel{eq:treex}
\end{equation}
It is clear that $\calf_n$ is the set of forests with depth less or equal to $n$.
{}From this observation, we see that
 $\calf_n$ form a linear ordered direct system: $\calf_n\supseteq \calf_{n-1}$,
and
\begin{equation}
 \calf =\cup_{n\geq 0} \calf_n=\dirlim \calf_n.
 \mlabel{lem:treex}
 \end{equation}
%Then we have $\lc \calf \rc\subseteq \calf.$
%

\subsection{Rota-Baxter operator on rooted forests}
\mlabel{ss:rbtree}
We note that $\bfk\,\calf$ with the product $\sqcup$ is also the free noncommutative nonunitary $\bfk$-algebra on the alphabet set $\calt$.
We are going to define, for each fixed $\lambda\in\bfk$,
another product $\shpr=\shpr_\lambda$ on $\bfk\,\calf$,
making it into a unitary Rota--Baxter algebra (of weight $\lambda$).
To ease notation, we will suppress $\lambda$.

We define $\shpr$ by giving a set map
$$
  \shpr: \calf \times \calf \to \bfk\, \calf
$$
and then extending it bilinearly. For this, we use the depth filtration
$\calf=\cup_{n\geq 0} \calf_n$ in Eq.~(\mref{lem:treex}) and apply induction on $i+j$ to define
$$
  \shpr: \calf_i \times \calf_j \to \bfk\, \calf.
$$
When $i+j=0$, we have $\calf_i=\calf_j=\langle \onetree\rangle$.
With the notation in Eq.~(\mref{eq:sqpower}), we define
\begin{equation}
     \shpr: \calf_0\times \calf_0 \to \bfk\,\calf,\
     \onetree^{\sqcup m} \shpr \onetree^{\sqcup n} :=\onetree^{\sqcup (m+n-1)}.
     \mlabel{eq:shprt0}
\end{equation}
For given $k\geq 0$, suppose that
$
  \shpr: \calf_i \times \calf_j \to \bfk\,\calf
$
is defined for $i+j\leq k$. Consider forests $F, F'$ with
$\depth(F)+\depth(F')=k+1$.

First assume that $F$ and $F'$ are trees.
Note that a tree is either $\onetree$ or is of the form $\lc
\oF\rc$ for a forest $\oF$ of smaller depth. Thus we can
define
\begin{equation}
F \shpr F' =\left \{ \begin{array}{ll}
%    \onetree, & {\rm\ if\ } T_1=T_2=\onetree, \\
    F, & {\rm\ if\ } F'=\onetree,\\
    F', & {\rm\ if\ } F=\onetree,\\
    \lc \lc \oF \rc \shpr \oF'\rc +\lc \oF\shpr \lc \oF' \rc \rc
            +\lambda\lc \oF \shpr \oF'\rc,
    & {\rm\ if\ } F=\lc \oF\rc, F'=\lc \oF'\rc,
    \end{array} \right .
\mlabel{eq:shprt1}
\end{equation}
since for the three products on the right hand of the third equation, the sums \allowdisplaybreaks{
\begin{eqnarray}
\depth(\lc\oF \rc)+ \depth(\oF'), \quad
\depth(\oF)+\depth(\lc \oF'\rc),
\quad \depth(\oF)+ \depth(\oF') \mlabel{eq:shprt1-1}
\end{eqnarray}}
are all less than or equal to $k$. Note that in either case,
$F\shpr F'$ is a tree or a sum of trees.

Now consider arbitrary forests $F=T_1\sqcup\cdots\sqcup T_b$ and
$F'=T'_1\sqcup \cdots\sqcup T'_{b'}$ with $\depth(F)+\depth(F')=k+1$.
We then define
\begin{equation}
F \shpr F'= T_1\sqcup\cdots \sqcup T_{b-1}\,
\sqcup\,(T_b\shpr T'_1)\,\sqcup \,
    T'_{2}\,\cdots\,\sqcup T_{b'}
\mlabel{eq:shprt2}
\end{equation}
where $T_b\shpr T'_1$ is defined by Eq.~(\mref{eq:shprt1}). By the remark after Eq.~(\mref{eq:shprt1-1}), $F\shpr F'$ is in $\bfk\,\calf$. This completes the definition of the set map $\shpr$
on $\calf \times \calf$.

As an example, we have
\begin{equation}
\td31 \shpr \tb2 = \lc \ta1\sqcup \ta1\rc \shpr \lc \ta1 \rc
= \lc (\ta1 \sqcup \ta1) \shpr \lc \ta1\rc\rc
+ \lc \lc \ta1 \sqcup \ta1\rc \shpr \ta1\rc
+ \lambda \lc (\ta1\sqcup\ta1) \shpr \ta1\rc
= \tg42 + \thj44 + \lambda \td31.
\mlabel{eq:treeex}
\end{equation}

We record the following simple properties of $\shpr$ for later applications.
\begin{lemma} Let $F,F',F''$ be forests.
\begin{enumerate}
\item
$(F \sqcup F')\shpr F'' =F\sqcup (F' \shpr F''), \quad
F''\shpr (F\sqcup F') =(F'' \shpr F)\sqcup F'.$
\mlabel{it:twoass}
\item
$ \leaf(F\shpr F')=\leaf(F)+\leaf(F')-1.$
\mlabel{it:leaf}
\end{enumerate}
\mlabel{lem:matcht}
\end{lemma}
So $\bfk\,\calf$ with the operations $\sqcup$ and $\shpr$ forms a
2-associative algebra in the sense of \cite{L-R3,Pi}.
\begin{proof}
(\mref{it:twoass}). Let $F=T_1\sqcup\cdots\sqcup T_b, $
$F'=T'_1\sqcup\cdots\sqcup T'_{b'}$ and
$F''=T''_1\sqcup\cdots\sqcup T''_{b''}$ be the decomposition of the forests into trees.
Since $\sqcup$ is an associative product, by Eq.~(\mref{eq:shprt2})
we have,
\allowdisplaybreaks{
\begin{eqnarray*}
 (F\sqcup F')\shpr F''&=&(T_1\sqcup\cdots\sqcup T_b\sqcup T'_1\sqcup \cdots
                            \sqcup T'_{b'}) \shpr (T''_1\sqcup T''_2\sqcup \cdots\sqcup  T''_{b''})\\
                      &=&T_1\sqcup\cdots\sqcup T_b\sqcup T'_1\sqcup \cdots \sqcup T'_{b'-1}
                            \sqcup (T'_{b'} \shpr T''_1)\sqcup T''_2\sqcup \cdots\sqcup  T''_{b''}\\
                      &=&(T_1\sqcup\cdots\sqcup T_b)\sqcup (T'_1\sqcup \cdots \sqcup T'_{b'-1}
                            \sqcup (T'_{b'} \shpr T''_1)\sqcup T''_2\sqcup \cdots\sqcup  T''_{b''})\\
                      &=&F\sqcup (F'\shpr F'').
\end{eqnarray*}}
The proof of the second equation is the same.

\noindent
(\mref{it:leaf}). We prove by induction on the sum
$m:=\depth(F)+\depth(F')$. When $m=0$, it follows from Eq.~(\mref{eq:shprt0}). Assume that the equation holds for all $F$ and $F'$ with $m\leq k$ and consider $F$ and $F'$ with $\depth(F)+\depth(F')=k+1$. If $F$ and $F'$ are trees, then the equation holds by Eq.~(\mref{eq:shprt1}), the induction hypothesis and the fact that $\leaf(\lc \oF\rc)=\leaf(\oF)$ for a forest $\oF$. Then for forests $F$ and $F'$, the equation follows from Eq.~(\mref{eq:shprt2}) and Eq.~(\mref{eq:leaf})
\end{proof}

Extending $\shpr$ bilinearly, we obtain
a binary operation
$$ \shpr: \bfk\,\calf \otimes \bfk\,\calf \to \bfk\, \calf.$$
For $F\in \calf$, we use the grafting operation to define
\begin{equation}
P_\calf(F)=\lc F\rc.
\mlabel{eq:RBt}
\end{equation}
Then $P_\calf$ extends to a linear operator on $\bfk\,\calf$.

The following is our first main result and will be proved in the next subsection.
\begin{theorem}
\begin{enumerate}
\item
    The pair $(\bfk\, \calf,\shpr)$ is a unitary associative algebra.
    \mlabel{it:algt}
\item
    The triple $(\bfk\,\calf,\shpr,P_\calf)$ is a unitary Rota--Baxter algebra of weight $\lambda$.
    \mlabel{it:RBt}
\end{enumerate}
\mlabel{thm:freet}
\end{theorem}

We next construct a nonunitary sub-Rota--Baxter algebra in $\bfk\,\calf$.

Let $\calfo$ be the subset of $\calf$ consisting of forests that are not $\onetree$ and do not contain any subtree $\lc \onetree \rc =\ \tb2\ $.
For example,
$$
 \td31  \quad\;
 \thj44 \quad\;
 \th43 \;\quad
$$
are in $\calfo$ while
$$
 \ta1 \;\quad
 \tb2   \quad\;
 \tc3   \quad\;
 \te4   \quad\;
 \tf41  \quad\;
 \tg42 \;\quad
$$
are not in $\calfo$.
Forests in $\calfo$ will be called the {\bf ladder-free forests}.

\begin{prop}
The submodule
$\bfk\,\calfo$ of $\bfk\,\calf$ is a nonunitary Rota--Baxter
subalgebra of $\bfk\,\calf$ under the product $\shpr$.
\mlabel{pp:freet0}
\end{prop}

\begin{proof}
We only need to check that $\bfk \calfo$ is closed under $\shpr$ and
$P_\calf=\lc\ \rc$. The following
lemma shows that $\bfk\, \calfo$ is closed under the Rota--Baxter
operator $P_\calf$.

\begin{lemma}
If $F$ is in $\calfo$, then $\lc F\rc$ does not contain $\lc
\onetree\rc$ and hence is in $\calfo$.
\mlabel{lem:oneb}
\end{lemma}
\begin{proof}
Let $F$ be in $\calfo$. Then $F$ does not contain $\lc\onetree\rc$.
In other words, none of the brackets $\lc B\rc$ in $F$ is of
the form $\lc \onetree\rc$. The only other brackets in $\lc F\rc$
is $\lc F\rc$ itself. So suppose $\lc F\rc$ contains a $\lc \onetree\rc$, then
we must have $\lc F\rc=\lc \onetree\rc$, implying $F=\onetree$.
This is a contradiction. So we have $\lc F\rc\in \calfo$.
\end{proof}

To prove that $\bfk\,\calfo$ is closed under the multiplication $\shpr$,
consider $F$ and $F'$ in $\calfo$. Since none of $F$ or $F'$ is
$\onetree$, we have $F\shpr F'\neq \onetree$. So the following lemma
completes the proof of Proposition~\mref{pp:freet0}.
\end{proof}

\begin{lemma}
If $F$ and $F'$ are in $\calfo$, then $F\shpr F'$ is either a forest that does not contain $\lc \onetree\rc$ or is a linear combination of forests that do not contain $\lc \onetree\rc$.
\mlabel{lem:twob}
\end{lemma}
\begin{proof}
Let $F=T_1\sqcup \cdots \sqcup T_b$ and
 $F'=T'_1\sqcup \cdots \sqcup T'_{b'}$.
We will prove the lemma using induction on $n:=\depth(T_b)+\depth(T'_1)$.

When $n=0$, we have $T_b=T'_1=\onetree$. Since none of $F$ or $F'$
is $\onetree$, we have $b>1$ and $b'>1$. So by Eq. (\mref{eq:shprt1}),
$$ F\shpr F'=T_1\sqcup \cdots \sqcup T_{b-1} \sqcup
\onetree \sqcup T'_2\sqcup \cdots \sqcup T'_{b'}.$$
Since neither $F$ nor $F'$ contains $\lc \onetree\rc$, none of
$T_i$ or $T'_j$ contains $\lc \onetree\rc$. Then none of the trees
on the right hand side contains $\lc \onetree \rc$.
So the right hand side does not contain $\lc \onetree\rc$,
as needed.

Let $k\geq 0$. Assume that the claim has been proved for $n\leq k$ and
let $F$ and $F'$ be in $\calfo$ with $n=k+1$. Then $n\geq 1$.
So at least one of $\depth(T_b)$
and $\depth(T'_1)$ is not zero. If one of them is zero, then the same
argument as in the $n=0$ case works using the first two cases of Eq.~(\mref{eq:shprt1}). If none of them is zero, then by the third case of Eq.~(\mref{eq:shprt1}), we have
$T_b = \lc \oF_b \rc$, $T'_{1}=\lc \oF'_1 \rc$ and
$$
  T_b \shpr T'_{b'}=\lc \lc \oF_b \rc \shpr \oF'_1 \rc
    + \lc \oF_b \shpr \lc \oF'_1 \rc\rc
    +\lambda \lc \oF_b \shpr \oF'_1 \rc.
$$
Since $T_b$ does not contain $\lc \onetree\rc$,
$\oF_b$ is not $\onetree$
and does not contain $\lc \onetree \rc$. So $\oF_b$ is in $\calfo$.
Similarly, $\oF'_1$ is in $\calfo$.
By the induction hypothesis, none of the terms
$\lc \oF_b \rc \shpr \oF'_1, \
 \oF_b \shpr \lc \oF'_1 \rc, \
 \oF_b \shpr \oF'_1
$
contains $\lc \onetree\rc$. Thus they are in $\bfk\,\calfo$.
By Lemma~\mref{lem:oneb},
the terms on the right hand side themselves do not contain
$\lc \onetree\rc$.
Therefore $T_b\shpr T'_{1}$ is a linear combination of terms
that do not contain $\lc \onetree\rc$.
Since $F$ and $F'$ do not contain $\lc\onetree\rc$, none of
$T_i$ and $T'_j$ contains $\lc \onetree\rc$.
By Eq.~(\mref{eq:shprt2}), we have
$$
F \shpr F'= T_1\sqcup\cdots \sqcup T_{b-1}\,
\sqcup\,(T_b\shpr T'_1)\,\sqcup \,
    T'_{2}\,\cdots\,\sqcup T_{b'}.
$$
Then $F\shpr F'$ is a linear combination of terms that do not
contain $\lc \onetree\rc$. This completes the induction.
\end{proof}

\subsection{The proof of Theorem~\mref{thm:freet}}
\mlabel{ss:prooft}
\begin{proof}
(\mref{it:algt}).
By Definition (\mref{eq:shprt1}), $\onetree$ is the identity under
the product $\shpr$.
So we just need to verify the associativity. For this we only need to verify
\begin{equation}
 (F\shpr F')\shpr F'' =F\shpr(F' \shpr F'')
 \mlabel{eq:asst}
\end{equation}
for forests $F,F',F''\in \calf$.
We will accomplish this by induction on the sum of the depths
$n:=\depth(F)+\depth(F')+\depth(F'')$. If $n=0$, then
all of $F,F',F''$ have depth zero and so are
in $\calf_0=\sqmon{\onetree}$, the sub-semigroup of $\calf$ generated by
$\onetree$. Then we have $F=\onetree^{\sqcup i}$,
$F'=\onetree^{\sqcup i'}$ and
$F''=\onetree^{\sqcup i''}$, for $i,i',i''\geq 1$.
Then the associativity follows from Eq.~(\mref{eq:shprt0}) since both
sides of Eq.~(\mref{eq:asst}) is $\onetree^{\sqcup (i+i'+i''-2)}$ in this case.

Let $k\geq 0$. Assume Eq.~(\mref{eq:asst}) holds for $n\leq k$ and assume that
$F,F',F''\in \calf$ satisfy
$n=\depth(F)+\depth(F')+\depth(F'')=k+1.$
We next reduce the breadths of the forests.

\begin{lemma}
If the associativity
$$
  (F \shpr F')\shpr F''= F\shpr (F' \shpr F'')
$$
holds when $F, F'$ and $F''$ are trees, then
it holds when they are forests.
\mlabel{lem:ellt}
\end{lemma}

\begin{proof}
We use induction on the sum of breadths $m:=b(F)+b(F')+b(F'')$.
Then $m\geq 3$. The case when $m=3$ is the assumption of the
lemma. Assume the associativity holds for $3\leq m \leq j$ and
take $F, F',F''\in \calf$ with $m = j+1.$ Then $j+1\geq 4$. So
at least one of $F,F',F''$ has breadth greater
than or equal to 2.

First assume $b(F)\geq 2$. Then $F=F_1\sqcup F_2$
with $F_1,\, F_2\in \calf$.
Thus by Lemma~\mref{lem:matcht},
$$
  (F\shpr F') \shpr F'' = ((F_1\sqcup F_2)\shpr F')\shpr F''
                           = (F_1 \sqcup (F_2 \shpr F'))\shpr F''
                           = F_1 \sqcup ((F_2 \shpr F') \shpr F'').
$$
Similarly,
$$
  F\shpr (F' \shpr F'') = (F_1\sqcup F_2)\shpr (F'\shpr F'')\\
                           = F_1 \sqcup (F_2 \shpr (F'\shpr F'')).
$$
Thus
$$
 (F\shpr F') \shpr F'' = F\shpr (F' \shpr F'')
$$
whenever
$$
  (F_2 \shpr F') \shpr F'' = F_2 \shpr (F'\shpr F'')
$$
which follows from the induction hypothesis.
A similar proof works if $b(F'')\geq 2$.

Finally if $b(F')\geq 2$, then $F'=F'_1\sqcup F'_2$
with $F'_1,\,F'_2\in \calf$.
Using Lemma~\mref{lem:matcht} repeatedly, we
have
$$ (F \shpr F')\shpr F''= (F \shpr (F'_1\sqcup F'_2)) \shpr F''
                         = ((F \shpr F'_1)\sqcup F'_2)\shpr F''
                         = (F\shpr F'_1)\sqcup (F'_2 \shpr F'').
$$
In the same way, we have
$
  F\shpr (F' \shpr F'')=(F\shpr F'_1)\sqcup (F'_2 \shpr F'').
$
This again proves the associativity.
\end{proof}

To summarize, our proof of the associativity (\mref{eq:asst}) has been reduced to the special case when the forests $F,F',F'' \in \calf$ are
chosen such that
\begin{enumerate}
\item
    $n:= \depth(F)+\depth(F')+\depth(F'')=k+1\geq 1$ with the assumption that the associativity holds when
    $n\leq k$, and
    \mlabel{it:sp1t}
\item
    the forests are of breadth one, that is, they are trees.
    \mlabel{it:sp2t}
\end{enumerate}

If either one of the trees is $\onetree$, the identity under
the product $\shpr$, then the associativity is clear.
So it remains to consider the case when $F,F',F''$ are all in
$\lc \calf \rc$. Then $F=\lc \oF\rc, F'=\lc \oF' \rc,
F''=\lc \oF''\rc$ with $\oF,\oF',\oF''\in \calf$. To deal with
this case, we prove the following
general fact on Rota--Baxter operators on not necessarily associative
algebras.

\begin{lemma}
Let $R$ be a $\bfk$-module with a multiplication $\cdot $ that is not
necessarily associative. Let $\lc\ \rc_R: R\to R$ be a $\bfk$-linear map
such that the Rota--Baxter identity holds:
\begin{equation}
\lc x\rc_R \cdot  \lc x'\rc_R
=\big\lc x \cdot  \lc x'\rc_R\big\rc_R
+\big\lc \lc x\rc_R \cdot  x'\big\rc_R
+\lambda \lc x \cdot  x'\rc_R,\
\forall\, x,\,x'\,\in R.
\mlabel{eq:rb2}
\end{equation}
%Define another multiplication $\odot_R$, called the {\bf double
%product} on $R$ by
%\begin{equation}
%x\odot_P\, y=x\odot P(y)+P(x)\odot y +\lambda x\odot y.
%\mlabel{eq:double}
%\end{equation}
%So we have $P(x\odot y)=P(x)\odot P(y)$.
Let $x,x'$ and $x''$ be in $R$. If
$$(x \cdot  x') \cdot  x''= x\cdot (x'\cdot x''),
$$
then we say that $(x,x',x'')$ is an {\bf associative triple}
for the product $\cdot$.
For any $y,y',y''\in R$, if all the triples
\begin{align}
&(y,y',y''),\ (\lc y\rc_R,y',y''),\ (y,\lc y'\rc_R,y''),\
(y,y',\lc y''\rc_R),\ (\lc y\rc_R,y',\lc y''\rc_R),
\mlabel{eq:triple1} \\
&(\lc y\rc_R, \lc y'\rc_R ,y''),\ (y,\lc y'\rc_R ,\lc y''\rc_R)
\mlabel{eq:triple2}
\end{align}
are associative triples for $\cdot$,
then $(\lc y\rc_R, \lc y'\rc_R, \lc y''\rc_R)$ is an associative
triple for $\cdot$.
\mlabel{lem:rbass}
\end{lemma}
\begin{proof}
Using
Eq.~(\mref{eq:rb2}) and bilinearity of the product $\cdot $, we
have \allowdisplaybreaks{
\begin{eqnarray*}
\lefteqn{(\lc y\rc_R \cdot  \lc y'\rc_R)\cdot \lc y''\rc_R
= \big( \lc \lc  y\rc_R \cdot   y '\rc_R + \lc y\cdot  \lc y'\rc_R\rc_R
    +\lambda \lc y\cdot   y' \rc \big ) \cdot  \lc  y''\rc_R }\\
&=& \lc\lc  y\rc_R \cdot   y'\rc_R \cdot  \lc y''\rc_R
    + \lc y\cdot  \lc  y'\rc_R \rc_R\cdot  \lc  y''\rc_R
    +\lambda \lc  y\cdot   y'\rc_R \cdot  \lc y''\rc_R \\
&=&  \lc\lc\lc  y\rc_R\cdot   y'\rc_R\cdot   y'' \rc_R
    + \lc\big(\lc y\rc_R \cdot   y'\big) \cdot  \lc y''\rc_R\rc_R
    +\lambda \lc\big(\lc y\rc_R \cdot  y'\big)\cdot   y''\rc_R\\
&& + \lc\lc y\cdot \lc y'\rc_R\rc_R \cdot   y''\rc_R
    + \lc\big( y\cdot \lc  y'\rc_R\big) \cdot \lc  y''\rc_R\rc_R
    +\lambda \lc\big( y\cdot  \lc  y'\rc_R \big) \cdot   y''\rc_R \\
&& + \lambda \lc \lc  y\cdot   y'\rc_R\cdot   y''\rc_R
    +\lambda \lc \big( y\cdot   y'\big)\cdot  \lc  y''\rc_R \rc_R
    + \lambda^2 \lc \big( y\cdot   y'\big) \cdot   y''\rc_R.
\end{eqnarray*}}
Applying the associativity of the second triple in
Eq.~(\mref{eq:triple2}) to $\big ( y\cdot \lc  y'\rc_R\big) \cdot \lc  y''\rc_R$ in the fifth term above and then using
Eq.~(\mref{eq:rb2}) again, we have
\allowdisplaybreaks{
\begin{eqnarray*}
&&(\lc y\rc_R \cdot  \lc y'\rc_R)\cdot \lc y''\rc_R\\
&=& \lc\lc\lc  y\rc_R\cdot  y'\rc_R\cdot
 y'' \rc_R
    + \lc\big(\lc y\rc_R \cdot  y'\big) \cdot \lc y''\rc_R\rc_R
    +\lambda \lc\big(\lc y\rc_R \cdot y'\big)\cdot  y''\rc_R\\
&& + \lc\lc y\cdot\lc y'\rc_R\rc_R \cdot  y''\rc_R
    + \lc y\cdot\lc\lc y'\rc_R\cdot y''\rc_R\rc_R
    + \lc y \cdot \lc y'\cdot \lc  y''\rc_R\rc_R\rc_R \\
&&  +\lambda \lc  y \cdot \lc y'\cdot  y''\rc_R\rc_R
    +\lambda \lc\big( y\cdot \lc  y'\rc_R \big) \cdot  y''\rc_R \\
&& + \lambda \lc \lc  y\cdot  y'\rc_R\cdot  y''\rc_R
    +\lambda \lc \big( y\cdot  y'\big)\cdot \lc  y''\rc_R \rc_R
    + \lambda^2 \lc \big( y\cdot  y'\big) \cdot  y''\rc_R.
\end{eqnarray*}}
By a similar calculation, we have
\allowdisplaybreaks{
\begin{eqnarray*}
\lefteqn{\lc y\rc_R \cdot \big(\lc y'\rc_R \cdot \lc y'' \rc_R\big)
= \lc\lc\lc y\rc_R\cdot  y'\rc_R\cdot  y''\rc_R
    + \lc \lc  y\cdot \lc y'\rc_R\rc_R \cdot  y''\rc_R}\\
&&  + \lambda \lc\lc y\cdot y'\rc_R\cdot y''\rc_R
  +\lc  y\cdot \lc \lc  y'\rc_R \cdot  y''\rc_R\rc_R
    + \lambda \lc  y\cdot \big(\lc y'\rc_R\cdot  y''\big)\rc_R\\
&&  + \lc\lc  y\rc_R\cdot \big( y'\cdot \lc  y''\rc_R \big) \rc_R
    + \lc  y \cdot \lc  y' \cdot \lc  y''\rc_R\rc_R\rc_R
    + \lambda \lc  y \cdot \big( y'\cdot \lc  y''\rc_R \big)\rc_R\\
&&  + \lambda\lc\lc y\rc_R\cdot \big(  y'\cdot  y''\big) \rc_R
    + \lambda \lc  y\cdot \lc  y'\cdot  y''\rc_R\rc_R
    + \lambda^2 \lc  y \cdot \big(  y'\cdot y''\big) \rc_R.
\end{eqnarray*}}
Now by the associativity of the triples in Eq.~(\mref{eq:triple1}),
the $i$-th term in the expansion of
$(\lc y\rc_R \cdot  \lc y'\rc_R)\cdot \lc y''\rc_R$
matches with the $\sigma(i)$-th
term  in the expansion of
$\lc y\rc_R \cdot \big(\lc y'\rc_R \cdot \lc y'' \rc_R\big)$.
Here the permutation $\sigma\in \Sigma_{11}$ is
\begin{equation}
\left ( \begin{array}{c} i\\\sigma(i)\end{array}\right)
= \left ( \begin{array}{ccccccccccc} 1&2&3&4&5&6&7&8&9&10&11\\
    1&6&9&2&4&7&10&5&3&8&11\end{array} \right ).
\mlabel{eq:sigmaa}
\end{equation}
This proves the lemma.
\end{proof}

To continue the proof of Theorem~\mref{thm:freet},
we apply Lemma~\mref{lem:rbass} to the situation where
$R$ is $\bfk\,\calf$ with the multiplication $\cdot = \shpr$,
the Rota--Baxter operator $\lc\ \rc_R=\lc\ \rc$ and the triple
$(y,y',y'')=(\oF,\oF',\oF'')$. By the induction hypothesis on $n$,
all the triples in Eq.~(\mref{eq:triple1}) and (\mref{eq:triple2}) are associative for $\shpr$.
So by Lemma~\mref{lem:rbass}, the triple $(F,F',F'')$ is associative for $\shpr$. This completes the induction and therefore
the proof of the first part of Theorem~\mref{thm:freet}.

(\mref{it:RBt}). We just need to prove that $P_\calf(F)=\lc F\rc$
is a Rota--Baxter operator of weight $\lambda$. This is immediate
from Eq. (\mref{eq:shprt1}).
\end{proof}

\section{Free Rota--Baxter algebras on a module or a set}
\mlabel{sec:freem}
We will construct the free unitary Rota--Baxter algebra on a $\bfk$-module or on a set by expressing elements in the Rota--Baxter algebra in terms of forests from Section~\mref{sec:rbt}, in addition with angles decorated by elements from the $\bfk$-module or set.
These decorated forests will be introduced in Section~\mref{ss:adecm}.
The free unitary Rota--Baxter algebra will be constructed
in Section~\mref{ss:freem}.
In Section~\mref{sec:freemo},
we also give a similar construction of free nonunitary Rota--Baxter
algebra on a $\bfk$-module in terms of the ladder-free forests
introduced in Proposition~\mref{pp:freet0}.
When the $\bfk$-module is taken to be free on a set,
we obtain the free unitary Rota--Baxter algebra on the set.
This will be discussed in Section~\mref{sec:freex}.

\subsection{Rooted forests with angular decoration by a module}
\mlabel{ss:adecm}
Let $M$ be a non-zero $\bfk$-module.
Let $F$ be in $\calf$ with $\leaf$ leafs.
We let $\tpow{F}{M}$ denote the tensor power $M^{\ot (\leaf-1)}$ labeled by $F$. In other words,
\begin{equation}
\tpow{F}{M}=\{ (F;\frakm)\ |\ \frakm\in M^{\ot (\leaf-1)}\}
\mlabel{eq:tpower}
\end{equation}
with the $\bfk$-module structure coming from the second component
and with the convention that $M^{\ot 0}=\bfk$.
We can think of $\tpow{F}{M}$ as the tensor power of $M$ with exponent $F$ with the usual tensor power $M^{\ot n}, n\geq 0$, corresponding to $\tpow{F}{M}$ when $F$ is the forest $\onetree ^{\sqcup (n+1)}$.

\begin{defn}
We call $\tpow{F}{M}$ the {\bf module of
the forest $F$ with angular decoration by $M$}, and call
$(F;\frakm)$, for $\frakm\in M^{\ot (\leaf(F)-1)}$, an {\bf angularly decorated forest} $F$ with the decoration tensor $\frakm$.
\mlabel{de:dectree}
\end{defn}
Also define the depth and breadth of $(F;\frakm)$ by
$$ \depth(F;\frakm)=\depth(F),\quad b(F;\frakm)=b(F).$$

Definition~\mref{de:dectree} is justified by the following tree interpretation of $\tpow{F}{M}$.
Let $(F;\frakm)$ be an angularly decorated forest with a pure tensor
$\frakm=a_1\ot \cdots \ot a_{\leaf-1}\in M^{\ot (\leaf-1)}$, $\leaf\geq 2$.
We picture $(F;\frakm)$ as the forest $F$ with its angles between
adjacent leafs (either from the same tree or from adjacent trees)
decorated by $a_1,\cdots,a_{\leaf-1}$ from the left most angle to
the right most angle.
If $\leaf(F)=1$, so $F$ is a ladder tree with only one leaf,
then $(F;a)$, $a\in\bfk$, is interpreted as
the multiple $a F$ of the ladder tree $F$.
For example, we have
$$
    \big( {\scalebox{1.15}{\tg42}}\ ;\ x     \big) = \begin{array}{l}\\[-.7cm] \xldec41r \end{array}, \quad
    \big( {\scalebox{1.15}{\thII43}}\ ;\ x\ot y \big) = \begin{array}{l}\\[-.5cm] \xyldec43 \end{array}, \quad
    \big( \ta1\ \sqcup\ {\scalebox{1.1}{\td31}}\ ;\ x\ot y \big)= \ta1 \sqcup_{x} \begin{array}{l}\\[-.3cm]
    \yldec31 \end{array}, \quad \big(\ta1\,;\ a \big)= a\, \ta1\,.
$$

When $\frakm=\sum_i \frakm_i$ is not a pure tensor,
but a sum of pure tensors $\frakm_i$ in $M^{\ot (\leaf-1)}$,
we can picture $(F;\frakm)$ as a sum $\sum_i  (F;\frakm_i)$ of
the forest $F$ with decorations from the pure tensors.
Likewise,
if $F$ is a linear combination $\sum_{i} c_i F_i$ of forests $F_i$ with
the same number of leaves $\leaf$ and if
$\frakm=a_1\ot \cdots\ot a_{\leaf-1}\in M^{\ot (\leaf-1)}$,
we also use $(F; \frakm)$ to denote the linear combination
$\sum_i c_i(F_i; \frakm)$.
For example,
$$
\big( {\scalebox{1.15}{\thII43}}+\, \ta1\, \sqcup\, {\scalebox{1.1}{\td31}}\ ;\ x\ot y \big)
= \begin{array}{l}\\[-.5cm] \xyldec43 \end{array} +\,  \ta1 \sqcup_{x} \begin{array}{l}\\[-.3cm]
    \yldec31 \end{array}
$$

Let $(F;\frakm)$ be an angular decoration of the forest $F$ by a pure tensor $\frakm$.
Let $F=T_1\sqcup \cdots \sqcup T_b$ be the decomposition of $F$
into trees. We consider the corresponding decomposition of
decorated forests. If $b=1$, then $F$ is a tree and $(F;\frakm)$
has no further decompositions. If $b>1$,
then there is the relation
$$ \leaf(F)=\leaf(T_1)+\cdots + \leaf(T_b).$$
Denote $\leaf_i=\leaf(T_i), 1\leq i\leq b$.
Then
$$(T_1;a_1\ot\cdots\ot a_{\leaf_1-1}),\
(T_2; a_{\leaf_1+1}\ot \cdots\ot a_{\leaf_1+\leaf_2-1}),
\cdots,
(T_b; a_{\leaf_1+\cdots+\leaf_{b-1}+1}\ot \cdots\ot a_{\leaf_1+\cdots+\leaf_b})
$$
are well-defined angularly decorated trees for the trees $T_i$ with $\leaf(T_i)>1$.
If $\leaf(T_i)=1$, then $a_{\leaf_{i-1}+\leaf_i-1}=a_{\leaf_{i-1}}$
and we use
the convention $(T_i;a_{\leaf_{i-1}+\leaf_i-1})=(T_i;\bfone)$.
With this convention, we have,
\begin{eqnarray*}
(F;a_1\ot\cdots\ot a_{\leaf-1})&=&
(T_1;a_1\ot \cdots\ot a_{\leaf_1-1})\sqcup_{a_{\leaf_1}}
(T_2; a_{\leaf_1+1}\ot\cdots\ot a_{\leaf_1+\leaf_2-1})
    \sqcup_{a_{\leaf_1+\leaf_2}}
\\
&&\cdots \sqcup_{a_{\leaf_1+\cdots+\leaf_{b-1}}}
(T_b; a_{\leaf_1+\cdots+\leaf_{b-1}+1}\ot \cdots\ot a_{\leaf_1+\cdots+\leaf_b}).
\end{eqnarray*}
We call this the {\bf standard decomposition} of $(F;\frakm)$ and
abbreviate it as
\begin{equation}
(F;\frakm)=(T_1;\frakm_1)\sqcup_{u_{1}}
(T_2;\frakm_2) \sqcup_{u_{2}} \cdots \sqcup_{u_{{b-1}}}
(T_b;\frakm_b).
\mlabel{eq:stdecm}
\end{equation}
In other words,
\begin{equation}
(T_i;\frakm_i)=\left \{ \begin{array}{ll}
(T_i;a_{\leaf_1+\cdots+\leaf_{i-1}+1}\ot \cdots \ot a_{\leaf_1+\cdots+\leaf_i-1}), & \leaf_i>1,i<b,\\
(T_i;a_{\leaf_1+\cdots+\leaf_{i-1}+1}\ot \cdots \ot a_{\leaf_1+\cdots+\leaf_i}), & \leaf_i>1,i=b,\\
(T_i;\bfone), & \leaf_i=1
\end{array} \right .
\end{equation}
and $u_i=a_{\leaf_1+\cdots+\leaf_i}.$
For example,
$$
    \big(\ta1\sqcup {\scalebox{1.15}{\tg42}} \sqcup {\scalebox{1.15}{\td31}}; v\ot x\ot w\ot y \big)
    = \big(\ta1;\bfone \big)\sqcup_v \big({\scalebox{1.15}{\tg42}};x)\sqcup_w \big({\scalebox{1.15}{\td31}};y \big)
    = \ta1 \sqcup_{v} \begin{array}{l}\\[-.7cm] \xldec41r \end{array} \sqcup_{w} \begin{array}{l}\\[-.3cm] \yldec31 \end{array}
$$
We display the following simple property for later applications.
\begin{lemma}
Let $F\neq \onetree$.
In the standard decomposition (\mref{eq:stdecm})
of $(F;\frakm)$, if $T_i=\onetree$ for some $1\leq i\leq b$, then $b>1$ and the corresponding factor $(T_i;\frakm_i)$ is $(T_i;\bfone)$.
\mlabel{lem:decone}
\end{lemma}
\begin{proof}
Let $F\neq \onetree$ and let $F=T_1\sqcup \cdots \sqcup T_b$ be its standard decomposition. Suppose $T_i=\onetree$ for some $1\leq i\leq b$ and $b=1$. Then $F=T_i=\onetree$, a contradiction. So $b>1$, and by our convention, $(T_i;\frakm_i)=(T_i;\onetree)$.
\end{proof}
\subsection{Free {Rota--Baxter algebra} on a module as decorated forests}
\mlabel{ss:freem}
We define the $\bfk$-module
$$ \ncsha(M)= \bigoplus_{F\in\, \calf} M^{\ot F}.$$
and define a product $\shprm$ on $\ncsha(M)$ by using
the product $\shpr$ on $\calf$ in Section~\mref{ss:rbtree}.

Let $T(M)=\oplus_{n\geq 0} M^{\ot n}$ be the tensor algebra and
let $\otm$ be its product, so for $\frakm\in M^{\ot n}$ and
$\frakm'\in M^{\ot n'}$, we have
\begin{equation}
\frakm \otm \frakm' =\left \{\begin{array}{ll}
 \frakm\ot\frakm'\in M^{\ot n+n'},& {\rm\ if\ } n>0,n'>0, \\
 \frakm \frakm'\in M^{\ot n'}, & {\rm if\ } n=0, n'>0,\\
\frakm'\frakm\in M^{\ot n}, & {\rm if\ } n>0, n'=0,\\
\frakm'\frakm\in \bfk, & {\rm if\ } n=n'=0.
\end{array} \right .
\mlabel{eq:multm}
\end{equation}
Here the products in the second and third case are scalar product
and in the fourth case is the product in $\bfk$.
In other words, $\otm$ identifies $\bfk\ot M$ and $M\ot \bfk$ with $M$ by the
structure maps $\bfk\ot M\to M$ and $M\ot \bfk \to M$ of the $\bfk$-module $M$.
\begin{defn}
\mlabel{de:shprm}
For tensors $D=(F;\frakm)\in \tpow{F}{M}$ and
$D'=(F';\frakm')\in \tpow{F'}{M}$, define
\begin{equation}
D\shprm D' = (F\shpr F'; \frakm\otm \frakm').
\mlabel{eq:shprm1}
\end{equation}
\end{defn}
The right hand side is well-defined since $\frakm\otm \frakm'$ has tensor degree
$$\deg(\frakm \otm \frakm')=\deg(\frakm)+\deg(\frakm')=\leaf(F)-1+\leaf(F')-1$$
which equals $\leaf(F\shpr F')-1$ by Lemma~\mref{lem:matcht}.(\mref{it:leaf}).
For example, from Eq.~(\mref{eq:treeex}) we have
$$
\begin{array}{l}\\[-.3cm] \xtd31 \end{array} \, \shprm \ \tb2
= \begin{array}{l}\\[-.7cm] \xldec41r \end{array}
+ \begin{array}{l}\\[-.7cm] \xthj44 \end{array}
+ \lambda \begin{array}{l}\\[-.3cm] \xtd31 \end{array}.
$$
By Eq. (\mref{eq:shprt0}) -- (\mref{eq:shprt2}), we have a more explicit expression.
\begin{equation}
 D \shprm D' = \left \{ \begin{array}{ll}
(\onetree; c c'),& {\rm\ if\ } D=(\onetree;c), D'=(\onetree;c'),\\
(F;c'\frakm), & {\rm\ if\ } D'=(\onetree, c'), F\neq \onetree,\\
(F'; c\frakm'), & {\rm\ if\ } D=(\onetree,c), F'\neq \onetree,\\
(F\shpr F'; \frakm \ot \frakm'),&  {\rm\ if\ } F\neq \onetree, F'\neq
\onetree. \end{array} \right .
\mlabel{eq:shprm}
\end{equation}

We can describe $\shprm$ even more explicitly in terms of the standard decompositions in Eq.~(\mref{eq:stdecm}) of $D=(F;\frakm)$ and $D'=(F';\frakm')$ for pure tensors $\frakm$
and $\frakm'$:
$$
D=(F;\frakm)=(T_1;\frakm_1)\sqcup_{u_{1}}
(T_2;\frakm_2) \sqcup_{u_{2}} \cdots \sqcup_{u_{{b-1}}}
(T_{b};\frakm_{b}),
$$
$$
D'=(F';\frakm')=(T'_1;\frakm'_1)\sqcup_{u'_{1}}
(T'_2;\frakm'_2) \sqcup_{u'_{2}} \cdots \sqcup_{u'_{{b'-1}}}
(T'_{b'};\frakm'_{b'}).
$$
Then by Eq.~(\mref{eq:shprt0}) -- (\mref{eq:shprt2}) and
Eq.~(\mref{eq:shprm1}) -- (\mref{eq:shprm}), it is easy to see that
the product $\shprm$ can be defined by induction on the sum of the depths $\depth=\depth(F)$ and $\depth'=\depth(F')$ as follows:
If $\depth+\depth'=0$, then $F=\onetree^{\sqcup i}$ and $F'=\onetree^{\sqcup j}$ for $i,j\geq 1$. If $i=1$, then $D=(F;\frakm)=(\onetree;c)=c(\onetree;\bfone)$ and we define
$ D\shprm D'=cD'=(F';c\frakm').$ Similarly define $D\shprm D'$ if $j=1$. If $i>1$ and $j>1$, then $(F;\frakm)=(\onetree;\bfone)\sqcup_{u_1}\cdots \sqcup_{u_{b-1}} (\onetree;\bfone)$ with $u_1,\cdots,u_{b-1}\in M$.
Similarly, $(F';\frakm')=(\onetree;\bfone)\sqcup_{u'_1}\cdots \sqcup_{u'_{b'-1}} (\onetree;\bfone).$ Then define
$$ (F;\frakm)\,\shprm \, (F';\frakm')=(\onetree;\bfone)\sqcup_{u_1}\cdots \sqcup_{u_{b-1}} (\onetree;\bfone) \sqcup_{u'_1} \cdots \sqcup_{u'_{b'-1}} (\onetree;\bfone).$$

Suppose $D\,\shprm\, D'$ has been defined for all $D=(F;\frakm)$ and $D'=(F';\frakm')$ with $\depth(F)+\depth(F')\leq k$ and consider
$D$ and $D'$ with $\depth(F)+\depth(F')=k+1$. Then we define
\begin{equation}
D \shprm D' =(T_1;\frakm_1)\sqcup_{u_{1}} \cdots \sqcup_{u_{{b-1}}}
\big((T_b;\frakm_b)\shprm (T'_1;\frakm'_1)\big) \sqcup_{u'_{1}}
 \cdots \sqcup_{u'_{{b'-1}}} (T'_{b'};\frakm'_{b'})
\mlabel{eq:shprm2}
\end{equation}
where
\begin{eqnarray}
&&(T_b;\frakm_b)\shprm (T'_1;\frakm'_1) \mlabel{eq:shprm3}\\
&=& \left \{ \begin{array}{ll}
(\onetree; \bfone), & {\rm if\ } T_b=T'_1=\onetree \ ({\rm so\ }
                \frakm_b=\frakm'_1=\bfone),\\
(T_b,\frakm_b), & {\rm if\ } T'_1=\onetree, T_b\neq \onetree, \\
(T'_1,\frakm'_1), & {\rm if\ } T'_1\neq\onetree, T_b= \onetree, \\
\lc (T_b;\frakm)\shprm (\oF'_1;\frakm')\rc + \lc (\oF_b;\frakm)\shprm (T'_1;\frakm') \rc  & \\
+\lambda \lc (\oF_b;\frakm)\shprm (\oF'_1;\frakm')\rc, &
    {\rm if\ } T'_1=\lc \oF'_1\rc \neq\onetree, T_b=\lc \oF_b\rc \neq \onetree.
\end{array} \right .
\notag
\end{eqnarray}
In the last case, we have applied the induction hypothesis on $\depth(F)+\depth(F')$ to define the terms in the brackets on the right hand side. Further, for $(F;\frakm)\in \tpow{F}{M}$, define $\lc (F;\frakm)\rc=(\lc F \rc;\frakm)$. This is well-defined since $\leaf(F)=\leaf (\lc F\rc)$.
\medskip

The product $\shprm$ is clearly bilinear.
So extending it biadditively, we obtain
a binary operation
$$
\shprm:  \ncsha (M)\otimes \ncsha(M) \to \ncsha(M).
$$
For $(F;\frakm) \in (F; M)$, define
\begin{equation}
 P_M(F; \frakm)=\lc (F;\frakm)\rc=(\lc F \rc\, ; \frakm)\in (\lc F\rc;M).
 \mlabel{eq:RBopm}
\end{equation}
As commented above, this is well-defined.
Thus $P_M$ defines a linear operator on $\ncsha(M)$. Note that the right
hand side is also
$(P_\calf(F); \frakm)$ with $P_\calf$ defined in
Eq.~(\mref{eq:RBt}). Let
\begin{equation}
  j_M: M \to \ncsha(M)
  \mlabel{eq:jm}
\end{equation}
be the $\bfk$-module map sending $a\in M$ to $(\onetree \sqcup \onetree;a)$.

\begin{theorem}
Let $M$ be a $\bfk$-module.
\begin{enumerate}
\item
    The pair $(\ncsha(M),\shprm)$ is a unitary associative algebra.
    \mlabel{it:algm}
\item
    The triple $(\ncsha(M),\shprm,P_M)$ is a unitary Rota--Baxter algebra of weight $\lambda$.
    \mlabel{it:RBm}
\item
    The quadruple $(\ncsha(M),\shprm,P_M,j_M)$ is the free unitary Rota--Baxter algebra of
    weight $\lambda$ on the module $M$. More precisely, for any unitary Rota--Baxter algebra $(R,P)$ and module morphism $f:M\to R$, there is a unique unitary Rota--Baxter algebra morphism
    $\free{f}: \ncsha(M) \to R$ such that $f=\free{f}\circ j_M.$
    \mlabel{it:freem}
\end{enumerate}
\mlabel{thm:freem}
\end{theorem}

\begin{proof}
(\mref{it:algm}) By definition, $(\onetree,\bfone)$ is the unit of
the multiplication $\shprm$.
For the associativity of $\shprm$ on $\ncsha(M)$
we only need to prove
$$
  (D \shprm D')\shprm D''=D\shprm (D' \shprm D'')
$$
for any angularly decorated forests $D=(F;\frakm) \in \tpow{F}{M},
D'=(F';\frakm')\in \tpow{F'}{M}$ and $D'' =(F'';\frakm'')\in \tpow{F''}{M}$.
Then by Eq.~(\mref{eq:shprm1}), we have
$$
(D\shprm D')\shprm D''=
\big( (F\shpr F')\shpr F''; (\frakm \otm \frakm') \otm \frakm''\big),$$
$$D\shprm (D'\shprm D'')=
\big( F\shpr (F'\shpr F''); \frakm \otm (\frakm' \otm \frakm'')\big).
$$
The first components of the two right
hand sides agree since the product $\shpr$ is associative by
Theorem~\mref{thm:freet}. The second component of the two right
hand sides agree because the product $\otm$ in Eq.~(\mref{eq:multm}) for the tensor algebra $T(M):=\bigoplus_{n\geq 0} M^{\ot n}$ is also associative. This proves the associativity of $\shprm$.
\medskip

(\mref{it:RBm}). The Rota--Baxter relation of $\lc\: \rc$ on
$\ncsha(M)$ follows from the Rota--Baxter relation of $\lc\: \rc$
on $\bfk\, \calf$ in Theorem~\mref{thm:freet}. More specifically,
it is the last equation in~(\mref{eq:shprm3}).
\medskip

(\mref{it:freem}). Let $(R,P)$ be a unitary Rota--Baxter algebra of weight
$\lambda$. Let $\ast$ be the multiplication in $R$ and let
$\bfone_R$ be its unit. Let $f: M\to R$ be a $\bfk$-module map. We will construct a
$\bfk$-linear map $\free{f}:\ncsha(M)\to R$ by defining
$\free{f}(D)$ for $D=(F;\frakm)\in \tpow{F}{M}$. We will achieve
this by induction on the depth $\depth(F)$ of $F$.

If $\depth(F)=0$, then $F=\onetree^{\sqcup i}$ for some $i\geq 1$.
If $i=1$, then $D=(\onetree; c)$, $c\in \bfk$.
Define $\free{f}(D)=c\bfone_R$.
In particular, define $\free{f}(\onetree;\bfone)=\bfone_R$.
Then $\free{f}$ sends the unit to the unit.
If $i\geq 2$, then $D=(F; \frakm)$ with $\frakm=a_1\ot \cdots\ot a_n
\in M^{\ot n}$ where $n+1$ is the number of leafs $\leaf(F)$.
Then we define
$\free{f}(\fraka)=f(a_1)\ast \cdots \ast f(a_n).$
In particular, $\free{f}\circ j_M = f$.

Assume that $\free{f}(D)$ has been defined for all $D=(F;\frakm)$ with $\depth(F)\leq k$ and let $D=(F;\frakm)$ with $\depth(F)=k+1$.
So $F\neq \onetree$. Let $D=(T_1;\frakm_1) \sqcup_{u_{1}} \cdots \sqcup_{u_{{b-1}}}
(T_b;\frakm_b)$ be the standard decomposition of $D$ given in
Eq. (\mref{eq:stdecm}). For each $1\leq i\leq b$, $T_i$ is a tree, so
it is either $\onetree$ or is of the form $\lc \oF_i\rc$ for another
forest $\oF_i$. By Lemma~\mref{lem:decone}, if $T_i=\onetree$, then $b>1$ and $\frakm_i=\bfone$. We accordingly define
\begin{equation}
\free{f}(T_i;\frakm_i)= \left\{ \begin{array}{ll}
\bfone_R, & {\rm if\ } T_i=\onetree,\\
P(\free{f}(\oF_i; \frakm_i)), & {\rm if\ } T_i=\lc \oF_i\rc.
\end{array} \right .
\mlabel{eq:homm0}
\end{equation}
In the later case, $(\oF_i;\frakm_i)$
is a well-defined angularly decorated forest since $\oF_i$ has
the same number of leafs as the number of leafs of $T_i$, and then
$\free{f}(\oF_i;\frakm_i)$ is defined
by the induction hypothesis since $\depth(\oF_i)=\depth(T_i)-1\leq k$.
Therefore we can define
\begin{equation}
\free{f}(D)= \free{f}(T_1;\frakm_1)\ast f(u_{1})\ast \cdots \ast
f(u_{{b-1}})\ast \free{f}(T_b;\frakm_b).
\mlabel{eq:homm}
\end{equation}

For any $D=(F;\frakm) \in \tpow{F}{M}$, we have
$P_M(D)=(\lc F\rc;\frakm) \in \ncsha(M)$, and
by the definition of $\free{f}$ in Eq.~(\mref{eq:homm0}) and  (\mref{eq:homm}), we have
\begin{equation}
 \free{f}(\lc D \rc)=P(\free{f}(D)).
\mlabel{eq:hom1-2m}
\end{equation}
So $\free{f}$ commutes with the Rota--Baxter operators.

Further, Eq. (\mref{eq:homm0}) and (\mref{eq:homm}) are clearly
the only way to define $\free{f}$ in order for $\free{f}$ to be a Rota--Baxter algebra homomorphism that extends $f$.

It remains to prove that the map $\free{f}$ defined in
Eq.~(\mref{eq:homm}) is indeed an algebra homomorphism. For this
we only need to check the multiplicativity
\begin{equation}
 \free{f} (D \shprm D)=\free{f}(D) \ast \free{f}(D')
 \mlabel{eq:hom2m}
\end{equation}
for all angularly decorated forests $D=(F;\frakm),D'=(F';\frakm')$ with pure tensors $\frakm$ and $\frakm'$.
Let
$$(F;\frakm)=(T_1;\frakm_1)\sqcup_{u_{1}}
(T_2;\frakm_2) \sqcup_{u_{2}} \cdots \sqcup_{u_{{b-1}}}
(T_b;\frakm_b)
$$
and
$$(F';\frakm')=(T'_1;\frakm'_1)\sqcup_{u'_{1}}
(T'_2;\frakm'_2) \sqcup_{u'_{2}} \cdots \sqcup_{u'_{{b'-1}}}
(T'_{b'};\frakm'_{b'})
$$
be their standard decompositions.

We first note that, since $\free{f}$ sends the identity
$(\onetree;\bfone)$ of $\ncsha(M)$ to the identity $\bfone_R$
of $R$, the multiplicativity is clear if either one of
$D$ or $D'$ is in $(\onetree;\bfk)$, that is, if either one of
$F$ or $F'$ is $\onetree$. So we only need to verify the
multiplicativity when $F\neq \onetree$ and $F'\neq \onetree$.

We further make the following reduction.
By Eq.~(\mref{eq:homm}) and Eq.~(\mref{eq:shprm2}), we have
\begin{eqnarray*}
\free{f}(D\shprm D') &=& \free{f}(T_1;\frakm_1) \ast f(u_{1})
\ast \cdots \ast f(u_{{b-1}}) \\
&& \ast
\free{f}\big((T_b;\frakm_b)\shprm (T'_1;\frakm'_1)\big) \ast f(u'_{1}) \ast
 \cdots \ast f(u'_{{b'-1}}) \ast
\free{f}(T'_{b'};\frakm'_{b'})
\end{eqnarray*}
and
\begin{eqnarray*}
\free{f}(D)\ast \free{f}(D') &=&
\free{f}(T_1;\frakm_1) \ast f(u_{1})
\ast \cdots \ast f(u_{{b-1}})\\
&& \ast \free{f}(T_b;\frakm_b)\ast \free{f}(T'_1;\frakm'_1)
\ast f(u'_{1}) \ast \cdots \ast f(u'_{{b'-1}}) \ast
\free{f}(T'_{b'};\frakm'_{b'}).
\end{eqnarray*}
We thus have
\begin{equation}
\free{f}((D;\frakm) \shprm (D';\frakm'))=
\free{f}(D;\frakm) \ast \free{f}(D';\frakm')
\mlabel{eq:mult1}
\end{equation}
if and only if
\begin{equation}
\free{f}((T_b,;\frakm_b)\shprm (T'_1;\frakm'_1))
=\free{f}(T_b;\frakm_b)\ast \free{f}(T'_1;\frakm'_1).
\mlabel{eq:mult2}
\end{equation}
So we only need to prove Eq.~(\mref{eq:mult2}).
For this we use induction on the sum of depths
$n:=\depth(T_b)+\depth(T'_1)$ of $T_b$ and $T'_1$. Then $n\geq 0$.
When $n=0$, we have $T_b=T'_1=\onetree$. So by Lemma~\mref{lem:decone}, we have $b>1, b'>1$, and
$$(T_b;\frakm_b)= (T'_1;\frakm'_1)
=(T_b;\frakm_b)\shprm (T'_1;\frakm'_1)=(\onetree;\bfone).$$
Then
$$ \free{f}(T_b;\frakm_b)= \free{f}(T'_1;\frakm'_1)
=\free{f}((T_b;\frakm_b)\shprm (T'_1;\frakm'_1))=\bfone_R.$$
Thus Eq.~(\mref{eq:mult2}) and hence Eq.~(\mref{eq:mult1}) holds.

Assume that the multiplicativity holds for
$D$ and $D'$
in $\tpow{\calf}{M}$ with $n=\depth(T_b)+\depth(T'_1)\leq k$ and take $D,D'\in \tpow{\calf}{M}$ with $n=k+1$. So $n\geq 1$. Then at least one of
$\depth(T_b)$ and $\depth(T'_1)$ is not zero. If exactly one of them is
zero, so exactly one of $T_b$ and $T'_1$ is $\onetree$,
then by Eq.~(\mref{eq:shprm3}),
$$ (T_b;\frakm_b)\shprm (T'_1;\frakm'_1) =\left \{\begin{array}{ll}
(T_b;\frakm_b), & {\rm\ if\ }  T'_1=\onetree, T_b\neq \onetree, \\
(T'_1;\frakm'_1), & {\rm\ if\ }  T'_1\neq\onetree, T_b= \onetree.
\end{array} \right .
$$
Then
$$
\free{f}((T_b;\frakm_b)\shprm (T'_1;\frakm'_1)) =\left\{ \begin{array}{ll}
\free{f}(T_b;\frakm_b), & {\rm\ if\ }  T'_1=\onetree, T_b\neq \onetree, \\
\free{f}(T'_1;\frakm'_1), & {\rm\ if\ }  T'_1\neq\onetree, T_b= \onetree.
\end{array} \right .
$$
Then Eq.~(\mref{eq:mult2}) and hence (\mref{eq:mult1})
holds since one factor in
$\free{f}(T_b;\frakm_b)\ast \free{f}(T'_1;\frakm'_1)$ is
$\bfone_R$.

If neither $\depth(T_b)$ nor $\depth(T'_1)$ is zero,
then $T_b=\lc \oF_b \rc$ and $T'_1=\lc \oF'_1 \rc$ for some forests
$\oF_b$ and $\oF'_1$ in $\calf$.
Then $(T_b;\frakm_b)=\lc (\oF_b;\frakm_b) \rc$
and  $(T'_1;\frakm'_1)=\lc (\oF'_1;\frakm'_1)\rc$.
We will take care of this case by the following lemma.

\begin{lemma}
Let $(R_1,P_1)$ and $(R_2,P_2)$ be not necessarily associative
$\bfk$-algebras $R_1$ and $R_2$ together with $\bfk$-linear
endomorphisms $P_1$ and $P_2$ that each satisfies the Rota--Baxter identity in
Eq.~(\mref{eq:RB}). Let $g:R_1\to R_2$ be a $\bfk$-linear map such that
\begin{equation}
g\circ P_1=P_2\circ g.
\mlabel{eq:rbinv}
\end{equation}
Let $x,y\in R_1$ be such that
\begin{equation}
g(xP_1(y))=g(x)\cdot g(P_1(y)),\ g(P_1(x)y)=g(P_1(x))\cdot g(y),\
g(xy)=g(x)\cdot g(y).
\mlabel{eq:rbhom1}
\end{equation}
Here we have suppressed the product in $R_1$ and denote the product
in $R_2$ by $\cdot$.
Then $g(P_1(x)P_1(y))=g(P_1(x))\cdot g(P_1(y)).$
\mlabel{lem:rbhom}
\end{lemma}
\begin{proof}
By the Rota--Baxter relations of $P_1$ and $P_2$,
Eq.~(\mref{eq:rbinv})
and Eq.~(\mref{eq:rbhom1}),
we have
{\allowdisplaybreaks
\begin{eqnarray*}
g(P_1(x)P_1(y))&=&g\big( P_1( P_1( x )   y )
    +P_1( x    P_1( y ) )
+\lambda P_1( x    y  )\big)\\
&=&g(P_1( P_1( x )   y ))
+ g(P_1( x    P_1( y ) ))
+g(\lambda P_1( x    y  ))\\
&=&P_2(g(P_1( x )   y ))
+ P_2(g(x    P_1( y ) ))
+ \lambda P_2(g(x    y  ))\\
&=&P_2(g(P_1( x ))\cdot  g(y ))
+ P_2(g(x ) \cdot  g( P_1( y ) ))
+ \lambda P_2(g(x ) \cdot  g(y ) )\\
&=&P_2(P_2(g(x ))\cdot  g(y ))
+ P_2(g(x ) \cdot  P_2(g(y )))
+ \lambda P_2(g(x ) \cdot  g(y ) )\\
&=& P_2(g(x ))\cdot  P_2(g(y ))\\
&=& g (P_1(x)) \cdot  g(P_1(y)).
\end{eqnarray*}
}
\end{proof}

Now we apply Lemma~\mref{lem:rbhom} to our proof with
$(R_1,P_1)=(\ncsha(M),\lc\ \rc)$, $(R_2,P_2)=(R,P)$ and
$g=\free{f}$. By the induction hypothesis, Eq.~(\mref{eq:rbhom1})
holds for $x=(\oF_b;\frakm_b)$ and $y=(\oF'_1;\frakm'_1)$.
Therefore by Lemma~\mref{lem:rbhom},
$\free{f}(T_b\shprm T'_1)=\free{f}(T_b)*\free{f}(T'_1)$.
Thus Eq.~(\mref{eq:mult1}) holds for $n=k+1$. This completes the
induction and the proof of Theorem~\mref{thm:freem}.
\end{proof}

%%%%%%%%%%%%%%%%%%%%%%%%%%%%%%%%%%%%%%%%%%%%%%%%%%%%%%%%%%%%%%%%%%%%%%%%%%%%
%%%%%%%%%%%%%%%%%%%%%%%%%%%%%%%%%%%%%%%%%%%%%%%%%%%%%%%%%%%%%%%%%%%%%%%%%%%%

\subsection{Free nonunitary Rota--Baxter algebra on a module}
\mlabel{sec:freemo} We now modify the construction of free unitary
Rota--Baxter algebras in Section~\mref{ss:freem} to obtain free
nonunitary Rota--Baxter algebras. Since the constructions are
quite similar, we will be brief for most parts except for the differences.

As in Proposition~\mref{pp:freet0}, we let $\calfo$ be the subset of
$\calf\backslash \{\onetree\}$ consisting of forests that do not contain any $\lc \onetree\rc =\ \tb2\ $.
For any $\bfk$-module $M$, define the $\bfk$-submodule
$$\ncshao(M) =\bigoplus_{F\in \calfo} \tpow{F}{M}$$
of $\ncsha(M)$.
We define a product $\shprm$ on $\ncshao(M)$ to be the restriction of $\shprm$ on $\ncsha(M)$. This product is well-defined since for $D=(F;\frakm)$ and $D'=(F';\frakm)$ with $F,F'\in \calfo$,
$F\shpr F'$ is in $\bfk\,\calfo$ by Proposition~\mref{pp:freet0}.
Thus by Eq.~(\mref{eq:shprm1}),
$D\shprm D' =(F \shpr F'; \frakm\otm \frakm')$
is in $\ncshao(M)$.

Also define $\lc\: \rc: \ncshao(M) \to \ncshao(M)$ to be the
restriction of $\lc\: \rc$ on $\ncsha(M)$. This again is well-defined
since by Proposition~\mref{pp:freet0}, $\lc \calfo\rc \subseteq \calfo$.
Then adapting the notation
and proof of Theorem~\mref{thm:freem}, we obtain

\begin{theorem}
Let $M$ be a $\bfk$-module.
\begin{enumerate}
 \item
    The pair $(\ncshao(M),\shprm)$ is a nonunitary associative algebra.
    \mlabel{it:algmo}
 \item
    The triple $(\ncshao(M),\shprm,P_M)$ is a nonunitary Rota--Baxter algebra of weight $\lambda$.
    \mlabel{it:RBmo}
 \item
    The quadruple $(\ncshao(M),\shprm,P_M,j_M)$ is the free nonunitary Rota--Baxter algebra of weight
    $\lambda$ on the $\bfk$-module $M$.
    \mlabel{it:freemo}
\end{enumerate}
\mlabel{thm:freemo}
\end{theorem}
\begin{proof}
(\mref{it:algmo}) and (\mref{it:RBmo}) are clear from (\mref{it:algm}) and (\mref{it:RBmo}) of Theorem~\mref{thm:freem}.

Part (\mref{it:freemo}) is proved in the same way as (\mref{it:freem}) of Theorem~\mref{thm:freem} with the following modification. Let $(R,\ast,P)$ be a nonunitary Rota--Baxter algebra. In the recursive definition of $\free{f}$ in Eq.~(\mref{eq:homm}), when $(T_i;\frakm_i)=(\onetree;\bfone)$, simply delete the factor $\free{f}(T_i;\frakm_i)$ instead of letting it be $\bfone_R$ which is not defined. Alternatively, augment $R$ to a unitary $\bfk$-algebra $\uni{R}=\bfk \bfone_R\oplus R$ with unit $\bfone_R$. Of course $\uni{R}$ can not be expected to be a Rota--Baxter algebra. But it does not matter since we only need the algebra structure on $\uni{R}$ to obtain a Rota--Baxter algebra structure on $R$.
For $D=(F;\frakm)\in \tpow{F}{M}$ with $F\in \calfo$, just define $\free{f}(D)$ as in Eq.~(\mref{eq:homm}).  Note that
$F$ has at least two leafs, so $\frakm$ is in $M^{\ot r}$ with
$r\geq 1$. Then it follows by induction that $\free{f}(D)$ is
always in $R$. Then the rest of the proof goes through.
\end{proof}

\subsection{Free Rota--Baxter algebra on a set}
\mlabel{sec:freex}
\mlabel{ss:freext}

Here we use the tree construction of free Rota--Baxter algebra on a module above to obtain a similar
construction of a free Rota--Baxter algebra on a set and display a canonical basis
of the free Rota--Baxter algebra in terms of forests decorated by
the set.

\begin{remark}{\rm
Either by the general principle of forgetful functors or by an easy direct check, the free Rota--Baxter algebra on a set $X$
is the free Rota--Baxter algebra on the
free $\bfk$-module
$M =\bfk\,X.$
Thus we can easily obtain a construction
of the free Rota--Baxter algebra on $X$ by decorated forests from the construction of
$\ncsha(M)$ in Section~\mref{ss:freem}.
}
\mlabel{rk:set}
\end{remark}

For any $n\geq 1$, the tensor power $M^{\ot n}$ has a natural basis
$X^n=\{ (x_1,\cdots,x_n)\ |\ x_i \in X, \ 1\leq i\leq n\}.$
Accordingly, for any rooted forest $F\in \calf$, with $\leaf=\leaf(F)\geq 2$, the set
$$X^F:=\{(F;(x_1,\cdots,x_{\leaf-1})):=(F;x_1\ot \cdots \ot x_{\leaf-1})\ |\ x_i\in X,\ 1\leq i\leq \leaf-1\}$$
form a basis of $\tpow{F}{M}$ defined in Eq.~(\mref{eq:tpower}).
Note that when $\leaf(F)=1$, $\tpow{F}{M}= \bfk\ F$ has a basis
$X^F:=\{(F;\bfone)\}$.
In summary, every $\tpow{F}{M}, F\in \calf,$ has a basis
\begin{equation}
X^{F}:=\{ (F; \vec{x})\ |\ \vec{x}\in X^{\leaf(F)-1}\},
\mlabel{eq:tpowerx}
\end{equation}
with the convention that $X^0=\{\bfone\}$. Thus the disjoint union
\begin{equation}
{X}^\calf:= \coprod_{F\in \calf} X^F.
\end{equation}
forms a basis of $$\ncsha(X):=\ncsha(M).$$
We call ${X}^\calf$ the set of {\bf angularly decorated rooted forests with
decoration set $X$}. As in Section~\mref{ss:adecm}, they can be
pictured as rooted forests with adjacent leafs decorated by
elements from $X$.

Likewise, for $(F;\vec{x})\in {X}^\calf$, the decomposition
(\mref{eq:stdecm}) gives the {\bf standard decomposition}
\begin{equation}
(F;\vec{x})=(T_1;\vec{x}_1)\sqcup_{u_{1}}
(T_2;\vec{x}_2) \sqcup_{u_{2}} \cdots \sqcup_{u_{{b-1}}}
(T_b;\vec{x}_b)
\mlabel{eq:stdecx}
\end{equation}
where $F=T_1\sqcup \cdots \sqcup T_b$ is the decomposition of
$F$ into trees and $\vec{x}$ is the vector concatenation of
the elements of
$\vec{x}_1, u_1,\vec{x}_2,\cdots,u_{b-1},\vec{x}_b$
which are not the unit $\bfone$.
As a corollary of Theorem~\mref{thm:freem}, we have
\begin{theorem}
For $D=(F;(x_1,\cdots,x_b))$, $D'=(F';(x'_1,\cdots,x'_{b'}))$
in $X^\calf$, define
\begin{equation}
D\shprm D'=\left \{\begin{array}{ll}
    (\onetree; \bfone), & {\rm\ if\ } F=F'=\onetree, \\
    D, & {\rm\ if\ } F'=\onetree, F\neq \onetree, \\
    D', & {\rm\ if\ } F=\onetree, F'\neq \onetree,\\
    (F\shpr F'; (x_1,\cdots,x_b,x'_1,\cdots,x'_{b'})), & {\rm \ if\ }
        F\neq \onetree, F'\neq \onetree,
        \end{array} \right.
\mlabel{eq:shprx}
\end{equation}
where $\shpr$ is defined in Eq.~(\mref{eq:shprt1}) and (\mref{eq:shprt2}).
Define
$$ P_X: \ncsha(X) \to \ncsha(X),\quad
P_X(F;(x_1,\cdots,x_b))=(\lc F\rc; (x_1,\cdots,x_b)),
$$
%,\quad
%(F;(x_1,\cdots,x_b))\in (F;X),$$
and
$$ j_X: X\to \ncsha(X),\quad
j_X(x)=(\onetree\sqcup \onetree; (x)),\quad x\in X.$$
Then the quadruple $(\ncsha(X),\shprm, P_X,j_X)$ is the free Rota--Baxter algebra
on $X$.
\mlabel{thm:freex}
\end{theorem}
\begin{proof}
The product $\shprm$ in Eq.~(\mref{eq:shprx}) is defined to be the
restriction of the product $\shprm$ in Eq.~(\mref{eq:shprm}) to
$X^\calf$. Since $X^\calf$ is a basis of $\ncsha(X)$,
the two products coincide. So
$\ncsha(X)$ and $\ncsha(M)$ are the same as Rota-Baxter algebras.
Then as commented in Remark~\mref{rk:set},
$\ncsha(X)$ is the free Rota--Baxter algebra on $X$.
\end{proof}
As with Theorem~\mref{thm:freemo}, the same proof there also gives
\begin{theorem}
The subalgebra $\ncshao(X)$ of $\ncsha(X)$ generated by the $\bfk$-basis
$X^{\calfo}:=\cup_{F\in \calfo} X^F$, with the same product $\shprm$, Rota--Baxter operator
$P_X$ and set map $j_X$, is the free nonunitary Rota--Baxter algebra on $X$.
\mlabel{thm:freexo}
\end{theorem}

%%%%%%%%%%%%%%%%%%%%%%%%%%%%%%%%%%%%%%%%%%%%%%%%%%%%
%%%%%%%%%%%%%%%%%%%%%%%%%%%%%%%%%%%%%%%%%%%%%%%%%%%%%%
\section{Unitarization of Rota--Baxter algebras}
\mlabel{sec:uni}
For any nonunitary algebra $A$ (even if
$A$ does have an identity), define $\uni{A}:=\bfk \oplus A$
with component wise addition and with product defined by
$$ (a,x) (b,y)=(ab, ay+bx+xy).$$
As is well-known, the unitarization of $A$ is $\uni{A}$ together
with the natural embedding
$$u_A: A \to \uni{A},\ x\mapsto (0,x).$$

To generalize this process to Rota--Baxter algebras turns
out to be much more involved since, after formally adding a unit
$\bfone$ to a nonunitary Rota--Baxter algebra $(A,P)$, we also need to add its images under the Rota--Baxter operator $P$ and its iterations,
such as $P(\bfone)$ and $P(xP(\bfone))$. Then it is not clear in general how
these new elements should fit together to form a Rota--Baxter algebra, except possibly in special cases (see Proposition~\mref{pp:idemp} below).
We will start with the unitarization of free Rota--Baxter algebras
and then take care of the case of a general Rota--Baxter algebra
by regarding it as a quotient of a free Rota--Baxter algebra.
Let us first give the definition.
\begin{defn} {\rm
Let $(A,P)$ be a
nonunitary Rota--Baxter $\bfk$-algebra. A unitarization of $A$ is
a unitary Rota--Baxter algebra $(\uni{A},\uni{P})$ with a
nonunitary Rota-Baxter algebra homomorphism $u_A:A\to \uni{A}$ such that for
any unitary Rota--Baxter algebra $B$ and a homomorphism $f: A\to
B$ of nonunitary Rota--Baxter algebras, there is a unique
homomorphism $\uni{f}: \uni{A} \to B$ of unitary Rota-Baxter algebras such that $f=\uni{f} \circ u_A$.}
\mlabel{de:unig}
\end{defn}

\subsection{Unitarization of free Rota--Baxter algebras}

Let $X$ be a set. Let $\ncsha(X)$ and $\ncshao(X)$ be the free
unitary and nonunitary Rota--Baxter algebras in Theorem~\mref{thm:freex} and Theorem~\mref{thm:freexo}. Let $\uni{j}_X:
X\to \ncsha(X)$ and $j_X:X\to \ncshao(X)$ be the canonical
embeddings. Regarding $\ncsha(X)$ as a nonunitary Rota--Baxter
algebra, then by the universal property of the free nonunitary
Rota--Baxter algebra $\ncshao(X)$, there is a unique homomorphism
$u_X:\ncshao(X) \to \ncsha(X)$ of nonunitary Rota--Baxter
algebras such that $\uni{j}_X=u_X\circ  j_X$.

\begin{theorem}
The unitary Rota--Baxter algebra
$\ncsha(X)$, with the homomorphism $u_X:\ncshao(X)\to \ncsha(X)$, is the unitarization
of the nonunitary Rota--Baxter algebra $\ncshao(X)$.
\mlabel{thm:unif}
\end{theorem}

\begin{proof}
Let $(B,Q)$ be a unitary Rota--Baxter algebra and let $f:
\ncshao(X)\to B$ be a homomorphism of nonunitary Rota--Baxter
algebras. Let $f'=f \circ j_X: X\to B$, then by
the freeness of the unitary Rota--Baxter algebra $\ncsha(X)$,
there is a unique homomorphism $\free{f}': \ncsha(X)\to B$ of
unitary Rota--Baxter algebras such that $f'=\free{f}'\circ \uni{j}_X$.

$$
      \xymatrix{ X \ar[rr]^{j_X} \ar[dd]^{\uni{j}_X} \ar[ddrr]^(.3){f'}
        && \ncshao(X) \ar[dd]_{f}  \ar[ddll]_(.3){u_X}\\ %\ar@/^3pc/[dd]^{\free{g}}
        &&\\
      \ncsha(X) \ar[rr]^{\free{f}'} \ar@/_1pc/[rr]_g && B}
$$
We have
$$
    \free{f}'\circ u_X \circ j_X = \free{f}'\circ \uni{j}_X
                                 =f' = f\circ j_X.
$$
By the freeness of $\ncshao(X)$, we have $\free{f}'\circ u_X =f.$
Suppose there is another unitary Rota--Baxter algebra homomorphism
$ g:\ncsha(X) \to B$ such that $ g\circ u_X = f$. Then
$$
  g\circ \uni{j}_X = g\circ u_X\circ j_X=f\circ j_X=f'=\uni{f}'\circ \uni{j}_X.
$$
So $g=\uni{f}'$ by the universal property of the free unitary
Rota--Baxter algebra $\ncsha(X)$.
\end{proof}

%%%%%%%%%%%%%%%%%%%%%%%%%%%%%%%%%%%%%%%%%%%%%%%%%%%%%%%%%%%%%

\subsection{Unitarization of Rota--Baxter algebras}

We now construct the unitarization of any given nonunitary
Rota--Baxter algebra $A$. We use the following diagram to keep
track of the maps that we will introduced below.
\begin{equation}
\xymatrix{
J\ar[dd]_{\incl} && X\ar[ddll]_{j_X} \ar[ddrr]^{\uni{j}_X}
    \ar[ddddddll]^(.7)g \ar[ddddddrr]^(.7){\uni{g}}
    && \uni{J}\ar[dd]_{\incl}\\
&&&&\\
\ncshao(X) \ar[rrrr]^{u_X} \ar[ddrr]^(.7)h \ar[dddd]_{\free{g}} &&&&
\ncsha(X) \ar[ddll]_(.7){\uni{h}} \ar[dddd]_{\free{\uni{g}}} \\
&&&&\\
&& B && \\
&&&&\\
A\cong \ncshao(X)/J \ar[uurr]_{f} \ar[rrrr]^{u_A} &&&&
    \uni{A}= \ncsha(X)/\uni{J} \ar@/_1pc/[uull]_{\uni{f}}
    \ar@/^1pc/[uull]^{\uni{f}'}
    }
\mlabel{eq:unitg}
\end{equation}

Let $X$ be a generating set of $A$ as a nonunitary Rota--Baxter
algebra with $g:X\hookrightarrow A$ being the inclusion map. Let
$\ncshao(X)$ be the free nonunitary Rota--Baxter algebra on $X$
with the canonical embedding $j_X: X\to \ncshao(X)$. Then there is
a unique nonunitary Rota--Baxter algebra homomorphism $\free{g}:
\ncshao(X) \to A$ such that $g=\free{g}\circ j_X$. Since $X$ is a
generating set of $A$, $\free{g}$ is surjective. So $A\cong
\ncshao(X)/J$ where $J$ is the kernel of $\free{g}$ and is a
Rota--Baxter ideal of $\ncshao(X)$. Recall from
Theorem~\mref{thm:unif} that we have the unitarization $u_X:
\ncshao(X) \to \ncsha(X)$. Let $\uni{J}$ be the Rota--Baxter ideal
of $\ncsha(X)$ generated by $u_X(J)$, and define
$$ \uni{A}=\ncsha(X)/\uni{J}$$
with $\free{\uni{g}}: \ncsha(X)\to \uni{A}$ being the quotient
Rota--Baxter homomorphism. Let $\uni{g}=\free{\uni{g}} \circ
\uni{j}_X$. Then $\free{\uni{g}}: \ncsha(X)\to \uni{A}$ is the
unique unitary Rota--Baxter algebra homomorphism induced from the
set map $\uni{g}$. So the notation $\free{\uni{g}}$ is justified.

Now since $u_X(J) \subseteq \uni{J}$, we have
$(\free{\uni{g}}\circ u_X)(J)=0$. Thus $\ker (\free{\uni{g}}\circ u_X)
\supseteq J$. Therefore, there is a unique homomorphism
$$ u_A: A\cong \ncshao(X)/J \to \uni{A} = \ncsha(X)/\uni{J}$$
of nonunitary Rota--Baxter algebras such that
$$ u_A \circ \free{g} = \free{\uni{g}} \circ u_X.$$

\begin{theorem}
With the above notations, the nonunitary Rota--Baxter algebra
homomorphism
$$ u_A: A\to \uni{A}$$
gives the unitarization of $A$.
\mlabel{thm:unitg}
\end{theorem}
By the uniqueness of the Rota--Baxter algebra unitarization, for a different
choices of the generating set $X$ of $A$, the unitarization we
obtain are isomorphic.
\begin{proof}
Let $B$ be a unitary Rota--Baxter algebra and let let $f:A\to B$
be a nonunitary Rota--Baxter algebra homomorphism. Let $h=f \circ \free{g}$. By Theorem~\mref{thm:unif}, there is a unique unitary Rota--Baxter algebra homomorphism $\uni{h}: \ncsha(X)\to B$ such that $\uni{h}\circ u_X=h$. Then $$ \ker \uni{h} \supseteq u_X(\ker h)\supseteq u_X(\ker \free{g})=J.$$
Since $\uni{h}$ is a Rota--Baxter ideal of $\ncsha(X)$ and
$\uni{J}$ is the Rota--Baxter ideal of $\ncsha(X)$ generated by
$J$, we must have $\ker \uni{h}\supseteq \uni{J}$. Therefore,
there is a unique
$$ \uni{f}: \uni{A} \to B$$
such that $\uni{h}=\free{\uni{g}}\circ \uni{f}$. Now
$$ \uni{f}\circ u_A \circ \free{g}=\uni{f}\circ \free{\uni{g}} \circ u_X
= \uni{h}\circ u_X = h= f \circ \free{g}.$$
Since $\free{g}$ is surjective, we have $\uni{f}\circ u_A=f$. So the existence
of $\uni{f}$ in Definition~\mref{de:unig} is proved.

To prove the uniqueness of $\uni{f}$, suppose there is also a unitary Rota--Baxter algebra homomorphism $\uni{f}':\uni{A} \to B$ such that $ \uni{f}'\circ u_A =f$. Then we have \allowdisplaybreaks{
\begin{eqnarray*}
 \uni{f}'\circ \free{\uni{g}} \circ u_X
 & = & \uni{f}'\circ u_A \circ \free{g}
 = f \circ \free{g}
 = \uni{f}\circ u_A \circ \free{g}
 =  \uni{f}\circ \free{\uni{g}} \circ u_X
 = \uni{h} \circ u_X
 = h.
 \end{eqnarray*}}
So $\uni{f}'\circ \free{\uni{g}}: \ncshao(X) \to B$, as well as $\uni{h}$ is
the unitarization of $h:\ncshao(X)\to B$. By the uniqueness of
this unitarization, proved in Theorem~\mref{thm:unif}, we have
$$\uni{f}'\circ \free{\uni{g}}=\uni{h}=\uni{f}\circ \free{\uni{g}}.$$
Since $\free{\uni{g}}$ is surjective, we have
$\uni{f}'=\uni{f}$, as needed.
\end{proof}

\subsection{Unitarization with idempotent Rota--Baxter operators}
We end our discussion on unitariness of Rota--Baxter algebras with a simple case.
\begin{prop}
Let $(R,P)$ be a Rota--Baxter algebra of weight $\lambda$ such that $P^2=-\lambda P$. The unitarization $\uni{R}:=\bfk \bfone \oplus R$ of $R$ together with the extension of $P$ to $\tilde{P}:\uni{R} \to \uni{R}$,
$$
    \tilde{P}(m,a):=\big(-\lambda m,P(a)\big),\;\;\forall m \in \bfk,\ a \in R,
$$
forms a unitary Rota--Baxter $\bfk$-algebra of weight $\lambda$ such that $\tilde{P}^2=-\lambda \tilde{P}$.
\mlabel{pp:idemp}
\end{prop}
Other results on such Rota--Baxter operators can be found in~\cite{A-M} where they are called pseudo-idempotent.
\begin{proof}
We first show that $\tilde{P}: \uni{R} \to \uni{R}$ satisfies the
Rota--Baxter relation of weight $\lambda$
\begin{equation}
      \tilde{P}(m,a)\tilde{P}(n,b)
      =\tilde{P}\big( (m,a)\tilde{P}(n,b)\big) + \tilde{P}\big( \tilde{P}(m,a)(n,b)\big)
      +\lambda \tilde{P}\big((m,a)(n,b)\big)
\mlabel{eq:rbid}
\end{equation}
for $(m,a),(n,b) \in \uni{R}$. For the left hand side, we have
 \allowdisplaybreaks{
\begin{eqnarray*}
    \tilde{P}(m,a)\tilde{P}(n,b) &=& \big(-\lambda m,P(a)\big) \big(-\lambda n,P(b)\big)\\
    &=& \big(\lambda^2 mn,-\lambda mP(b) - \lambda nP(a) + P(a)P(b)\big)\nonumber\\
    &=& \big(\lambda^2 mn, -\lambda mP(b) - \lambda nP(a) + P(aP(b)) + P(P(a)b) + \lambda P(ab)\big).\nonumber
\end{eqnarray*}}
For the right hand side, we have
\allowdisplaybreaks{
\begin{eqnarray*}
    \tilde{P}\big( (m,a)\tilde{P}(n,b)\big)
%    &=& \tilde{P}\big( (m,a)(-\lambda n,P(b)) \big)\\
    &=& \tilde{P}(-\lambda mn,mP(b) -\lambda na + aP(b)) \nonumber\\
    &=& \big( \lambda^2 mn, mP^2(b) -\lambda nP(a) + P(aP(b))\big)\\
    &=& \big( \lambda^2 mn, -\lambda mP(b) -\lambda nP(a) + P(aP(b))\big),\nonumber
\end{eqnarray*}}
where we have used idempotency of $P$ in the second equality. For the
other terms we similarly find \allowdisplaybreaks{
\begin{eqnarray*}
\tilde{P}\big( \tilde{P}(m,a)(n,b)\big)
%    &=& \tilde{P}\big( (-\lambda m,P(a))(n,b) \big)\\
    &=& \big( \lambda^2 mn, -\lambda mP(b) -\lambda nP(a) + P(P(a)b)\big),\nonumber\\[0.3cm]
\tilde{P}\big((m,a)(n,b)\big)
%    &=& \tilde{P}\big( mn,na + mb + ab \big)\\
    &=& \big( -\lambda mn,mP(b) + nP(a) + P(ab) \big).\nonumber
\end{eqnarray*}}
{}From these equations, Eq.~(\mref{eq:rbid}) is immediately verified. %Note that the addition in $\tilde{R}=\bfk \oplus R$ is defined componentwise.

Finally,
$$ \tilde{P}^2(m,a)=\tilde{P}(-\lambda m; P(a))=((-\lambda)^2 m;P^2(a))=(\lambda^2 m; -\lambda P(a))=-\lambda \tilde{P}(m,a).$$
\end{proof}

%%%%%%%%%%%%%%%%%%%%%%%%%%%%%%%%%%%%%%%%%%%%%%%%%%%%%%%%%%%%%%%%%%%%%%%%%%%%%

%
%\addcontentsline{toc}{section}{\numberline {}References}
%

\end{document}